\documentclass{amsart}     
\usepackage{graphicx}
\usepackage{amsmath}
\usepackage{amssymb}
\usepackage{amsfonts}
\usepackage{latexsym}
\usepackage{oldgerm}

\newcommand{\smallmattwo}[4]{\left(\begin{smallmatrix} #1 & #2 \\ #3 &#4			       \end{smallmatrix}\right)}

\newcommand{\mattwo}[4]{\left(\begin{array}{cc} #1 & #2 \\ #3 &#4			       \end{array}\right)}

\begin{document}

\title{On the period  of the Ikeda lift for $U(m,m)$}

\author{Hidenori Katsurada \\
Muroran Institute of Technology 27-1 Mizumoto Muroran 050-8585, Japan
             }

\maketitle

\begin{abstract}
Let $K={\bf Q}(\sqrt{-D})$ be an imaginary quadratic field with discriminant $-D,$  and $\chi$ the Dirichlet character corresponding to the extension $K/{\bf Q}.$ 
Let $m=2n$ or $2n+1$ with $n$ a positive integer. Let $f$ be a primitive form of weight $2k+1$ and character $\chi$ for $\varGamma_0(D),$ 
  or a primitive form of weight $2k$ for $SL_2({\bf Z})$ according as $m=2n,$ or $m=2n+1.$ For such an $f$   let $I_m(f)$ be the lift of $f$ to the space of modular forms of weight $2k+2n$ and character $\det^{-k-n}$ for the Hermitian modular group $\varGamma_K^{(m)}$ constructed by Ikeda.
  We then express the period $\langle I_m(f), I_m(f) \rangle$ of $I_m(f)$ 
in terms of special values of the adjoint $L$-function of $f$ and its  twist by the character $\chi.$ 
This proves the conjecture concerning the period of the Hermitian Ikeda lift proposed by Ikeda.  

\end{abstract}

\section{Introduction} 
It is an important and interesting problem to consider the relation between the period of an elliptic modular form and that of its lift.
Here, we say that $F$ is a lift of an elliptic modular form $f$ if $F$ or the adelization of $F$ is a Hecke eigenform in the space of Siegel cusp forms or Hermitian cusp forms whose
certain $L$-function is expressed in terms of $L$-functions related to $f.$ 
There are several results concerning this problem in the Siegel modular form case (cf. \cite{B-D-S}, \cite{K-S}). 
This type of period relation sometimes gives rise to congruence between the lift and non-lift, and  are important also from the view point of arithmetic geometry (cf. \cite{B-D-S}, \cite{Br}, \cite{Kat1}). 
In \cite{K-K4}, we proved a conjecture on the period of the Duke-Imamoglu-Ikeda lift (DII lift) proposed by Ikeda \cite{Ik2}. As a result, in \cite{Kat2}, we characterized prime ideals giving congruence between the DII lift and non-DII lift. (See also 
\cite{Br-K}.)  Klosin \cite{Klo} gave the congruence between the Hermitian Maass lift and non-Hermitian Maass lift using the period relation 
in \cite{Ik3}.  In this paper we prove a result similar to \cite{K-K4} for the period of the lift of an elliptic modular form to the space of Hermitian modular forms constructed by Ikeda. This also proves Ikeda's conjecture in \cite{Ik3} with some modification. 

Let $K={\bf Q}(\sqrt{-D})$ be an imaginary quadratic field with discriminant $-D$, and  $\chi$ the Kronecker character corresponding to the extension $K/{\bf Q}.$ 
Let $k$ be a non-negative integer. Then for a primitive form $f \in {\mathfrak S}_{2k+1}(\varGamma_0(D),\chi)$  Ikeda \cite{Ik3} constructed a lift $I_{2n}(f)$ of $f$ to the space of modular forms of weight $2k+2n$ and a character $\det^{-k-n}$ for the Hermitian group $\varGamma_K^{(2n)}$ of degree $m.$ This is a generalization of the Maass lift considered by Kojima \cite{Koj}, Gritsenko \cite{G}, Krieg \cite{Kri}, Oda \cite{O}, and Sugano \cite{Su}. Similarly for a primitive form $f \in {\mathfrak S}_{2k}(SL_2({\bf Z}))$  he  constructed a lift $I_{2n+1}(f)$ of $f$ to the space of modular forms of weight $2k+2n$ and a character $\det^{-k-n}$ for $\varGamma_K^{(2n+1)}.$ For the rest of this section, let $m=2n$ or $m=2n+1.$ We then call $I_{m}(f)$  the Ikeda lift of $f$ for $U(m,m)$ or the Hermitian Ikeda lift of degree $m.$ Then our main result (Theorem 2.1) can be stated as follows:

{\it The period  $\langle I_m(f), I_m(f) \rangle$ of $I_m(f)$ is  expressed as 
$$L(1,f,{\rm Ad}) \prod_{i=2}^{m} L(i,f,{\rm Ad},\chi^{i-1})L(i,\chi^i)$$
 up to elementary factor, where $L(s,f,{\rm Ad},\chi^{i-1})$ is the "modified twist" of the adjoint $L$-function of $f$ by $\chi^{i-1},$ and $L(i,\chi^i)$ is the Dirichlet $L$-function for $\chi^i.$}
 
  This  result was already obtained in the case $m=2$, and  was conjectured  in general case by Ikeda \cite{Ik3}.  

We note that $I_m(f)$ is not likely to be a theta lift  except in the case $m=2,$  and therefore  the method in \cite{Ra} cannot be applied to prove our main result. The method we use is similar to that in the proof of the main result of \cite{K-K4} and to give  an explicit formula of  the  Dirichlet series of Rankin-Selberg type associated to $I_m(f)$, and to compare its residue with $\langle I_m(f), I_m(f) \rangle$. We explain it more precisely. In Section 3, we
 consider the  Dirichlet series $R(s,I_m(f))$  of Rankin  Selberg type associated with $I_m(f).$
  For the precise definition, see Section 3. This type of  Dirichlet series  was studied by Shimura \cite{Sh2}  for a classical Hermitian modular form $F$ of weight $2k+2n.$ In particular we can express its residue at $2k+2n$ in terms of the period of $F$ (cf. Proposition 3.1).  Thus to prove Theorem 2.1, we have to get an explicit formula of  $R(s,I_m(f))$  in terms of $L(s,f,{\rm Ad},\chi^i).$  To get it, in Section 4, we reduce our computation to  a computation of certain formal power series $\hat   H_{m,p}(d;X,Y,t)$ in  $t$ associated with local Siegel series similarly to \cite{K-K4}  (cf. Theorem 4.1). 

Section 5 is devoted to the computation of them. This computation is similar to that in \cite{K-K4}, but we should be careful in dealing with the case where $p$ is ramified in $K$. After such an elaborate computation, we can get explicit formulas of $\hat H_{m,p}(d;X,Y,t)$ for all prime numbers $p$ (cf. Theorem 5.5.4). In Section 6, by using explicit formulas for $\hat H_{m,p}(d;X,Y,t)$, we immediately get an explicit formula of $R(s,I_m(f))$ (cf. Theorems 6
.1 and 6.2) and by taking the residue of it at $2k+2n$ we prove the Theorem 2.1.  

We  note that we can give a similar period relation for the adelic Ikeda lift, and we can apply it to a problem concerning congruence between the adelic Ikeda lifts and Hecke eigenforms not coming from the adelic Ikeda lifts. These  will be discussed in  subsequent papers.

\bigskip
  
\qquad {\bf Notation.}  Let $R$ be a commutative ring. We denote by $R^{\times}$ and $R^*$  the semigroup of non-zero elements of $R$ and the unit group of $R,$  respectively. For a subset $S$ of $R$ we denote by $M_{mn}(S)$ the set of
$(m,n)$-matrices with entries in $S.$ In particular put $M_n(S)=M_{nn}(S).$ 
  Put $GL_m(R) = \{A \in M_m(R) \ | \ \det A \in R^* \},$ where $\det
A$ denotes the determinant of a square matrix $A$. Let $K_0$ be a field, and  $K$ a quadratic extension of $K_0,$ or $K=K_0 \oplus K_0.$ In the latter case, we regard $K_0$ as a subring of $K$ via the diagonal embedding.  We also identify $M_{mn}(K)$ with $M_{mn}(K_0) \oplus M_{mn}(K_0)$ in this case.
If $K$ is a quadratic extension of $K_0,$ let  $\rho$ be the non-trivial automorphism of $K$ over $K_0,$ and if $K=K_0 \oplus K_0,$ let $\rho$ be the automorphism of $K$ defined by $\rho(a,b)=(b,a)$ for $(a,b) \in K_0.$ We sometimes write $\overline {x}$ instead of $\rho(x)$ for $x \in K$ in both cases. Let $R$ be a subring of $K.$ For an $(m,n)$-matrix $X=(x_{ij})_{m \times n}$ write $\overline X=(\overline{x_{ij}})_{m \times n}$ and $X^*={}^t \overline X$,  and for an $(m,m)$-matrix
$A$, we write $A[X] = X^* A X.$ Let ${\mathrm{Her}}_n(R)$ denote
the set of Hermitian  matrices of degree $n$ with entries in
$R$, that is the subset of $M_n(R)$ consisting of matrices $X$ such that $X^*=X.$   Then a Hermitian matrix $A$ of degree $n$ with entries in $K$ is said to be semi-integral over $R$  if ${\rm tr}(AB) \in K_0 \cap R$ for any $B \in {\mathrm{Her}}_n(R),$ where ${\rm tr}$ denotes the trace of a matrix. We denote by $\widehat{{\mathrm{Her}}}_n(R)$ the set of semi-integral matrices of degree $n$ over $R.$ 

  For a subset $S$ of
$M_n(R)$ we denote by $S^{\times}$ the subset of $S$
consisting of non-degenerate matrices. If $S$ is a subset of ${\mathrm{Her}}_n({\bf C})$ with ${\bf C}$ the field of complex  numbers, we denote by $S^+$ the subset of $S$ consisting of positive definite matrices. The group 
  $GL_n(R)$ acts on the set ${\mathrm{Her}}_n(R)$ from the right in the following way:
  $$ GL_n(R) \times {\mathrm{Her}}_n(R) \ni (g,A) \longrightarrow g^* Ag \in {\mathrm{Her}}_n(R).$$
  Let $G$ be a subgroup of $GL_n(R).$ For a $G$-stable subset ${\mathcal B}$ of ${\mathrm{Her}}_n(R)$ we denote by ${\mathcal B}/G$ the set of equivalence classes of ${\mathcal B}$ under the action of $G.$ We sometimes identify ${\mathcal B}/G$ with a complete set of representatives of ${\mathcal B}/G.$ We abbreviate ${\mathcal B}/GL_n(R)$ as ${\mathcal B}/\sim$ if there is no fear of confusion. Two Hermitian matrices $A$ and $A'$ with
entries in $R$ are said to be $G$-equivalent and write $A \sim_{G} A'$ if there is
an element $X$ of $G$ such that $A'=A[X].$ 
For square matrices $X$ and $Y$ we write $X \bot Y = \mattwo{X}{O}{O}{Y}.$

We put ${\bf e}(x)=\exp(2 \pi \sqrt{-1} x)$ for $x \in {\bf C},$ and for a prime number $p$ we denote by ${\bf e}_p(*)$ the continuous additive character of ${\bf Q}_p$ such that ${\bf e}_p(x)= {\bf e}(x)$ for $x \in {\bf Z}[p^{-1}].$ 

For a prime number $p$ we denote by ${\rm ord}_p(*)$ the additive valuation of ${\bf Q}_p$ normalized so that ${\rm ord}_p(p)=1,$ and put  $|x|_p=p^{-{\rm ord}_p(x)}.$ Moreover we denote by $|x|_{\infty}$ the absolute value of $x \in {\bf C}.$

\section{Period of the Ikeda lift  for $U(m,m)$}
For a positive integer $N$ let $\varGamma_0(N)=\{ \mattwo{a}{b}{c}{d} \in SL_2({\bf Z}) \ | \ c \equiv 0 \ {\rm mod} \ N  \},$
and for a Dirichlet character $\psi$ mod $N,$ we denote by ${\mathfrak M}_{l}(\varGamma_0(N),\psi)$ the space of modular forms of weight $l$ for 
$\varGamma_0(N)$ and nebentype $\psi,$ and by ${\mathfrak S}_{l}(\varGamma_0(N),\psi)$ its subspace consisting of cusp forms. We simply write
 ${\mathfrak M}_{l}(\varGamma_0(N),\psi)$ (resp. ${\mathfrak S}_{l}(\varGamma_0(N),\psi)$)  as ${\mathfrak M}_{l}(\varGamma_0(N))$ (resp. as ${\mathfrak S}_{l}(\varGamma_0(N))$) if $\psi$ is the trivial character.  

Throughout the paper, we fix an imaginary quadratic extension $K$ of ${\bf Q}$ with the discriminant $-D,$ and denote by ${\mathcal O}$ the ring of integers in $K.$  
For a prime number $p$ put $K_p=K \otimes {\bf Q}_p,$ and ${\mathcal O}_p={\mathcal O} \otimes {\bf Z}_p.$ Then $K_p$ is a quadratic extension of ${\bf Q}_p$ or $K_p \cong {\bf Q}_p \oplus {\bf Q}_p.$  In the former case, for  $x \in K_p,$ we denote by $\overline {x}$ the conjugate of $x$ over ${\bf Q}_p.$ In the latter case, we identify $K_p$ with ${\bf Q}_p \oplus {\bf Q}_p,$ and for $x=(x_1,x_2) \in {\bf Q}_p \oplus {\bf Q}_p,$ we put $\overline {x}=(x_2,x_1).$ For $x \in K_p$ we define the norm $N_{K_p/{\bf Q}_p}(x)$ by $N_{K_p/{\bf Q}_p}(x)=x\overline {x},$ and 
 put $\nu_{K_p}(x)={\rm ord}_p(N_{K_p/{\bf Q}_p}(x)),$ and $|x|_{K_p}=|N_{K_p/{\bf Q}_p}(x)|_p.$ Moreover put $|x|_{K_{\infty}}=|x \overline{x}|_{\infty}$ for $x \in {\bf C}.$

For a non-degenerate Hermitian matrix or alternating matrix $T$ with entries in $K,$ let ${\mathcal U}_T$ be the unitary group defined over ${\bf Q},$ whose group ${\mathcal U}_T(R)$ of $R$-valued points is given by
$${\mathcal U}_T(R)=\{ g \in GL_{m}(R \otimes K) \ | \ {}^t\overline{g}Tg = T \}$$
for any ${\bf Q}$-algebra $R,$ where $g \mapsto \overline {g}$ denotes the automorphism of $M_n(R \otimes K)$ induced by the non-trivial automorphism of $K$ over ${\bf Q}.$ 
We also define the special unitary group $\mathcal {SU}_T$ over ${\bf Q}_p$ by $\mathcal {SU}_T={\mathcal U}_T \cap R_{K/{\bf Q}}(SL_m),$ where $R_{K/{\bf Q}}$ is the Weil restriction. In particular we write ${\mathcal U}_{J_m}$ as ${\mathcal U}^{(m)}$ or $U(m,m),$ where $J_m=\smallmattwo{O}{-1_m}{1_m}{O}.$ Then   
$${\mathcal U}^{(m)}({\bf Q})=\{M \in GL_{2m}(K)  \ | \  J_m[M]= J_m  \}.$$
 Put 
$$\varGamma^{(m)}=\varGamma_K^{(m)}={\mathcal U}^{(m)}({\bf Q}) \cap GL_{2m}({\mathcal O}). $$
  Let ${\mathfrak H}_m$ be the Hermitian upper half-space defined by 
$${\mathfrak H}_m=\{Z \in M_m({\bf C}) \ | \ {1 \over 2\sqrt{-1}} (Z-Z^*) \ {\rm is \ positive \ definite} \}.$$
The group ${\mathcal U}^{(m)}({\bf R})$ acts on ${\mathfrak H}_m$ by
$$ g \langle Z \rangle =(AZ+B)(CZ+D)^{-1} \ {\rm for} \ g=\smallmattwo{A}{B}{C}{D} \in {\mathcal U}^{(m)}({\bf R}), Z \in {\mathfrak H}_m.$$
We also put $j(g,Z)=\det (CZ+D)$ for such $Z$ and $g.$
  Let $l$ be an integer. For a subgroup $\varGamma$ of  ${\mathcal U}^{(m)}({\bf Q})$ which is commensurable with $\varGamma^{(m)}$ and a character $\psi$ of $\varGamma,$ we denote by ${\mathfrak M}_{l}(\varGamma,\psi)$  the space of holomorphic modular forms of weight $l$  with character $\psi$ for $\varGamma.$   We denote by ${\mathfrak S}_{l}(\varGamma,\psi)$ the subspace of ${\mathfrak M}_{l}(\varGamma,\psi)$ consisting of cusp forms. In particular, if $\psi$ is the character of $\varGamma$ defined  by $\psi(\gamma)=(\det \gamma)^{-l}$ for $\gamma \in \varGamma,$ we write ${\mathfrak M}_{2l}(\varGamma,\psi)$ as  ${\mathfrak M}_{2l}(\varGamma,\det^{-l}),$ and so on.  
Write the variable $Z$ on ${\mathfrak H}_m$ as $Z=X+\sqrt{-1}Y$ with $X,Y \in {\mathrm{Her}}_m({\bf C}).$ We can identify  ${\mathrm{Her}}_m({\bf C})$ with ${\bf R}^{m^2}$ through the map $X=(x_{ij}) \longrightarrow (x_{ii},{\rm Re}(x_{ij}), {\rm Im}(x_{ij}) \ (i < j)),$ and define a measure $dX$ on ${\mathrm{Her}}_m({\bf C})$ by pulling back the standard measure on ${\bf R}^{m^2}.$ Similarly we define a measure $dY$ on ${\mathrm{Her}}_m({\bf C})$ in the same way as above. 
  For two cusp forms $F$ and $G$ of weight $l$ with respect to $\varGamma^{(m)}$ with character $\chi$  we define the Petersson scalar product $\langle F,G \rangle$ by 
$$\langle F,G \rangle=\int_{\varGamma^{(m)} \backslash {\mathfrak H}_m} F(Z)\overline {G(Z)} (\det Y)^{l-2m} dXdY,$$
 where $X={Z+^t\overline{Z} \over 2},$ and $Y={Z-{}^t \overline{Z} \over 2 \sqrt{-1}}.$
 We call $\langle F, F \rangle$ the period of $F.$ 
 Similarly for two elements $f,g \in {\mathfrak S}_{l}(\varGamma_0(N), \psi),$
we define the Petersson scalar product $\langle f,g \rangle$ by 
$$\langle f,g \rangle=[SL_2({\bf Z}):\varGamma_0(N)]^{-1}\int_{\varGamma \backslash {\mathfrak H}} f(z)\overline {g(z)} y^{l-2} dxdy,$$
where ${\mathfrak H}$ is the complex upper half space. 

Now we consider adelic modular forms.  
Let ${\bf A}$ be the adele ring of ${\bf Q},$ and  ${\bf A}_{f}$ the non-archimedian factor of ${\bf A}.$ 
  Let $h=h_K$ be a class number of $K.$ 
Let $G^{(m)}={\rm Res}_{K/{\bf Q}}(GL_m),$ and $G^{(m)}({\bf A})$ be the adelization of $G^{(m)}.$ Moreover  put ${\mathcal C}^{(m)}=\prod_p GL_m({\mathcal O}_p).$ Let ${\mathcal U}^{(m)}({\bf A})$ be the adelization of ${\mathcal U}^{(m)}.$ We define the  compact subgroup $\mathcal K_0^{(m)}$ of ${\mathcal U}^{(m)}({\bf A}_f)$ by ${\mathcal U}^{(m)}({\bf A}) \cap \prod_p GL_{2m}({\mathcal O}_p),$ where $p$ runs over all rational primes. 
 Then we have
$${\mathcal U}^{(m)}({\bf A})=\bigsqcup_{i=1}^h {\mathcal U}^{(m)}({\bf Q}) \gamma_i \mathcal K_0 ^{(m)}{\mathcal U}^{(m)}({\bf R})$$
with some subset $\{\gamma_1,...,\gamma_h \}$ of ${\mathcal U}^{(m)}({\bf A}_f).$ We can take $\gamma_i $ as
$$\gamma_i =\mattwo{t_i}{0}{0}{t_i^{* -1}},$$ 
where $\{ t_i \}_{i=1}^h =\{(t_{i,p}) \}_{i=1}^h $ is a certain subset of   $G^{(m)}({\bf A}_f)$ such that $t_1=1,$ and
 $$G^{(m)}({\bf A})=\bigsqcup_{i=1}^h G^{(m)}({\bf Q}) t_i G^{(m)}({\bf R}){\mathcal C}^{(m)}.$$
  Put $\varGamma_i={\mathcal U}^{(m)}({\bf Q}) \cap \gamma_i{\mathcal K}_0 \gamma_i^{-1}{\mathcal U}^{(m)}({\bf R}).$ Then for an element $(F_1,...,F_h) \in \bigoplus_{i=1}^h  {\mathfrak M}_{2l}(\varGamma_i,\det^{-l}),$ we define $(F_1,...,F_h)^\sharp$ by
$$(F_1,...,F_h)^{\sharp}(g)=F_i(x\langle {\bf i} \rangle)j(x,{\bf i})^{-2l} (\det x)^l$$
for $g=u\gamma_ix \kappa$ with $u \in {\mathcal U}^{(m)}({\bf Q}),x \in {\mathcal U}^{(m)}({\bf R}),\kappa \in {\mathcal K}_0.$  We denote by ${\mathcal M}_l({\mathcal U}^{(m)}({\bf Q}) \backslash {\mathcal U}^{(m)}({\bf A}), \det^{-l})$ the space of automorphic forms obtained in this way. We also put
$${\mathcal S}_{2l}({\mathcal U}^{(m)}({\bf Q}) \backslash {\mathcal U}^{(m)}({\bf A}),{\det}^{-l})=\{ (F_1,...,F_h)^{\sharp} \ | \ F_i \in {\mathfrak S}_{2l}(\varGamma_i,{\det}^{-l}) \}.$$
We can define the Hecke operators which act on the space \\
${\mathcal M}_{2l}({\mathcal U}^{(m)}({\bf Q}) \backslash {\mathcal U}^{(m)}({\bf A}), \det^{-l}).$ For the precise definition of them, see \cite{Ik3}.

 Let $\widehat {\mathrm{Her}}_m({\mathcal O})$ be the set of semi-integral Hermitian matrices over ${\mathcal O}$ of degree $m$ as in the Notation. We note that $A \in {\mathrm{Her}}_m(K)$ belongs to  $\widehat {\mathrm{Her}}_m({\mathcal O})$ if and only if its diagonal components are rational integers and $\sqrt{-D} A \in M_m({\mathcal O}).$   

For a non-degenerate Hermitian matrix $B$ with entries in $K_p$ of degree $m,$ put
 $\gamma(B)=(-D)^{[m/2]}\det B.$  Let  
$\widehat {\mathrm{Her}}_m({\mathcal O}_p)$ be the set of semi-integral matrices over ${\mathcal O}_p$ of degree $m$ as in the Notation. We put $\xi_p=1,-1,$ or
 $0$ according as $K_p = {\bf Q}_p \oplus {\bf Q}_p, K_p$ is
 an unramified quadratic extension of ${\bf Q}_p,$ or $K_p$
 is a  ramified quadratic extension of ${\bf Q}_p.$ 
For $T \in \widehat {\mathrm{Her}}_{m}({\mathcal O}_p)^{\times}$ we define the local Siegel series $b_p(T,s)$ by   $$b_p(T,s)=\sum_{R \in {\mathrm{Her}}_n(K_p)/{\mathrm{Her}}_n({\mathcal O}_p)} {\bf e}_p({\rm tr}(TR))p^{-{\rm ord}_p(\mu_p(R))s},$$ where $\mu_p(R)=[R{\mathcal O}_p^m+{\mathcal O}_p^m:{\mathcal O}_p^m]^{1/2}.$  

\noindent
{\bf Remark.}  In \cite{Kat3}, we defined $\mu_p(R)$ as $\mu_p(R)=[R{\mathcal O}_p^m+{\mathcal O}_p^m:{\mathcal O}_p^m].$  However,
it should be defined as above.

We remark that there exists a unique polynomial 
 $F_p(T,X)$ in $X$ such that 
 $$b_p(T,s)=F_p(T,p^{-s})\prod_{i=0}^{[(m-1)/2]}(1-p^{2i-s})\prod_{i=1}^{[m/2]} (1-\xi_p p^{2i-1-s}) $$ 
(cf. Shimura \cite{Sh1}). 
We then define a Laurent polynomial $\widetilde F_p(T,X)$ as
$$\widetilde F_p(T,X)=X^{-{\rm ord}_p(\gamma(T))}F_p(T,p^{-m}X^{2}).$$
We remark that we have 
$$\widetilde F_p(T,X^{-1})=(-D,\gamma(T))_p \widetilde F_p(T,X)  \qquad {\rm if } \ m \ {\rm is \ even},$$ 
$$\widetilde F_p(T,\xi_pX^{-1})=\widetilde F_p(T,X)  \qquad {\rm if } \ m \ {\rm is \ even \ and} \ p \nmid D,$$ 
and
$$\widetilde F_p(T,X^{-1})=\widetilde F_p(T,X)  \qquad {\rm if } \ m \ {\rm is \ odd}$$ (cf. \cite{Ik3}). Here $(a,b)_p$ is the Hilbert symbol of $a,b \in {\bf Q}_p^{\times}.$ Hence we have
$$\widetilde F_p(T,X)=(-D,\gamma(B))_p^{m-1} X^{{\rm ord}_p(\gamma(T))}F_p(T,p^{-m}X^{-2}).$$
Now we put
 $${\widehat {\mathrm{Her}}}_{m}({\mathcal O})_i^+=\{T \in {\mathrm{Her}}_m(K)^+ \ | \  t_{i,p}^*T t_{i,p} \in \widehat {\mathrm{Her}}_m({\mathcal O}_p) \ {\rm for \ any } \ p \}.$$

 Let  $k$ be a non-negative integer. First let $m=2n$ be a positive even integer and let 
$$f(z)=\sum_{N=1}^{\infty}a(N){\bf e}(Nz)$$
 be a  primitive form in ${\mathfrak S}_{2k+1}(\varGamma_0(D),\chi).$  For a prime number $p$ not dividing $D$ let $\alpha_p \in {\bf C}$ such that $\alpha_p+\chi(p)\alpha_p^{-1}=p^{-k}a(p),$ and for $p \mid D$ put $\alpha_p=p^{-k}a(p).$  We note that $\alpha_p \not=0$ even if $p | D.$ Then for the Kronecker character $\chi$ we  define Hecke's $L$-function $L(s,f,\chi^i)$ twisted by $\chi^i$ as 
  $$L(s,f,\chi^i)=\prod_{p \nmid D}\{(1-\alpha_p p^{-s+k}\chi(p)^i)(1-\alpha_p^{-1} p^{-s+k} \chi(p)^{i+1})\}^{-1}$$
$$\times 
\left\{\begin{array}{ll}
\prod_{p \mid D} (1-\alpha_p p^{-s+k})^{-1} & {\rm if} \ i \ {\rm is \ even} \\
\prod_{p \mid D} (1-\alpha_p^{-1} p^{-s+k})^{-1} & {\rm if} \ i \ {\rm is \ odd}. \end{array} 
\right. $$
 In particular, if $i$ is even,  we sometimes write $L(s,f,\chi^i)$ as $L(s,f)$ as usual. Moreover  we define a Fourier series
 $$I_m(f)(Z)= \sum_{T \in \widehat {\mathrm{Her}}_{m}({\mathcal O})^+}a_{I_{m}(f)}(T){\bf e}({\rm tr}(TZ)),$$
 where
 $$a_{I_{2n}(f)}(T)= |\gamma (T)|^{k} \prod_p \widetilde F_p(T,\alpha_p^{-1}).$$
Next let $m=2n+1$ be a positive odd integer and let 
$$f(z)=\sum_{N=1}^{\infty}a(N){\bf e}(Nz)$$
 be a  primitive form in ${\mathfrak S}_{2k}(SL_2({\bf Z})).$  For a prime number $p$ let $\alpha_p \in {\bf C}$ such that $\alpha_p+\alpha_p^{-1}=p^{-k+1/2}a(p).$ Then we  define Hecke's $L$-function $L(s,f,\chi^i)$ twisted by $\chi^i$ as   $$L(s,f,\chi^i)$$
  $$=\prod_{p }\{(1-\alpha_p p^{-s+k-1/2}\chi(p)^{i})(1-\alpha_p^{-1} p^{-s+k-1/2} \chi(p)^{i})\}^{-1}.$$
 In particular, if $i$ is even we write $L(s,f,\chi^i)$ as $L(s,f)$ as usual. We define a Fourier series
 $$I_{2n+1}(f)(Z)= \sum_{T \in \widehat {\mathrm{Her}}_{2n+1}({\mathcal O})^+} a_{I_{2n+1}(f)}(T){\bf e}({\rm tr}(TZ)),$$
 where
 $$a_{I_{2n+1}(f)}(T)= {|\gamma(T)|}^{k-1/2} \prod_p \widetilde F_p(T,\alpha_p^{-1}).$$
 \noindent
 {\bf Remark.}  In \cite{Ik3}, Ikeda defined $\widetilde F_p(T,X)$ as 
$$\widetilde F_p(T,X)=X^{{\rm ord}_p(\gamma(T))}F_p(T,p^{-m}X^{-2}),$$
and we define it by replacing $X$ with $X^{-1}$ in this paper. This change does not affect the results.

 Then Ikeda \cite{Ik3} showed  the following: 
 
 \bigskip
  
{\it Let $m=2n$ or $2n+1.$ Let  $f$ be a primitive form in 
${\mathfrak S}_{2k+1}(\varGamma_0(D),\chi)$ or in ${\mathfrak S}_{2k}(SL_2({\bf Z}))$ according as $m=2n$ or $m=2n+1.$ 
Then $I_m(f)(Z)$ is an element of ${\mathfrak S}_{2k+2n}(\varGamma^{(m)},\det^{-k-n}).$ 
  \it}
 
\bigskip

To state our main result, put 
\[\Gamma_{{\bf R}}(s)=\pi^{-s/2}\Gamma(s/2)\]
and 
\[\Gamma_{{\bf C}}(s)=\Gamma_{{\bf R}}(s)\Gamma_{{\bf R}}(s+1).\]
We note that
\[\Gamma_{{\bf C}}(s)=2(2\pi)^{-s} \Gamma(s).\]
For an integer $i$ let  $L(s,\chi^i)=\zeta(s)$ or $L(s,\chi)$ according as $i$ is even or odd, where $\zeta(s)$ and $L(s,\chi)$
 are Riemann's zeta function, and Dirichlet $L$-function for $\chi,$ respectively, and put 
$$	\widetilde{\Lambda}(s,\chi^i) = \Gamma_{{\bf C}}(s) L(s,\chi^i).$$ 
For a primitive form $f$ in ${\mathfrak S}_{2k+1}(\varGamma_0(D),\chi),$
 we define the adjoint $L$-function $L(s,\, f,\, {\rm Ad})$  and its twist $L(s,\, f,\, {\rm Ad},\chi)$ by $\chi$ as 
$$L(s,f,{\rm Ad})=\prod_{p \nmid D} \{(1-\alpha_p^2 \chi(p)p^{-s})(1-\alpha_p^{-2}\chi(p) p^{-s})(1-p^{-s}) \}^{-1}\prod_{p \mid D}(1-p^{-s})^{-1},$$
and
$$L(s,f,{\rm Ad},\chi)=\prod_{p \nmid D} \{(1-\alpha_p^2 p^{-s})(1-\alpha_p^{-2}p^{-s})(1-\chi(p)p^{-s}) \}^{-1}$$ 
$$\times \prod_{p \mid D} \{(1-\alpha_p^2 p^{-s})(1-\alpha_p^{-2}p^{-s})\}^{-1}.$$ 
For a primitive form $f$ in ${\mathfrak S}_{2k}(SL_2({\bf Z})),$
 we define the adjoint $L$-function $L(s,\, f,\, {\rm Ad})$  and its twist $L(s,\, f,\, {\rm Ad},\chi)$ by $\chi$ as 
$$L(s,f,{\rm Ad})=\prod_{p } \{(1-\alpha_p^2 p^{-s})(1-\alpha_p^{-2}p^{-s})(1-p^{-s}) \}^{-1},$$
and
$$L(s,f,{\rm Ad},\chi)=\prod_{p} \{(1-\alpha_p^2 \chi(p) p^{-s})(1-\alpha_p^{-2}\chi(p)p^{-s})(1-\chi(p)p^{-s}) \}^{-1}.$$ 
Let $f$ be a primitive form in ${\mathfrak S}_{2k+1}(\varGamma_0(D),\chi)$ or in ${\mathfrak S}_{2k}(SL_2({\bf Z}))$ according as $m=2n$ or $m=2n+1.$ 
We then  put 
$$L(s,\, f,\, {\rm Ad},\chi^i)= \left\{\begin{array}{ll} L(s,\, f, {\rm Ad}) \ & \ {\rm if} \ i \ {\rm is \ even} \\
L(s,f,{\rm Ad},\chi) & \ {\rm if} \ i \ {\rm is \ odd}
\end{array}
\right.$$
Moreover put
$$\widetilde{\Lambda}(s,\, f,\, {\rm Ad},\chi^i)= \Gamma_{{\bf C}}(s)\Gamma_{{\bf C}}(s+l-1)L(s,\, f,\, {\rm Ad},\chi^i),$$
where $l=2k+1$ or $l=2k$ according as $f \in {\mathfrak S}_{2k+1}(\varGamma_0(D),\chi)$ or $f \in {\mathfrak S}_{2k}(SL_2({\bf Z})).$
 Let
 $Q_D$ be the set of prime divisors of $D.$ For each prime $q \in Q_D,$ put $D_q=q^{{\rm ord}_q(D)}.$ We define a Dirichlet character $\chi_q$ by
 $$\chi_q(a)= 
\begin{cases} \chi(a') &  \ {\rm if} \ (a,q)=1\\
 0 &  \ {\rm if} \ q|a
\end{cases} ,$$
 where $a'$ is an integer such that
 $$a' \equiv a \ {\rm mod} \ D_q \quad \ {\rm and} \ a' \equiv 1 \ {\rm mod} \ DD_q^{-1}.$$
  For a subset $Q$ of $Q_D$ put
  $\chi_Q=\prod_{q \in Q} \chi_q$ and $\chi'_Q=\prod_{q \in Q_D, q \not\in Q} \chi_q.$ Here we make the convention that $\chi_Q=1$ and $\chi'_Q=\chi$ if $Q$ is the empty set. Let 
  $$f(z)=\sum_{N=1}^{\infty} c_f(N){\bf e}(Nz)$$
   be a primitive form in ${\mathfrak S}_{2k+1}(\varGamma_0(D),\chi).$ Then there exists a primitive form 
   $$f_Q(z)=\sum_{N=1}^{\infty} c_{f_Q}(N){\bf e}(Nz)$$
   such that 
   $$c_{f_Q}(p)=\chi_Q(p)c_f(p) \ {\rm for} \ p \not\in Q$$
   and
   $$c_{f_Q}(p)=\chi'_Q(p)\overline{c_f(p)} \ {\rm for} \ p \in Q.$$

 Then our main result in this paper is:
 
 \bigskip

\noindent
{\bf Theorem 2.1.} {\it
{\rm (1)} Let $m=2n$ be a positive even integer. For a primitive form $f$ in ${\mathfrak S}_{2k+1}(\varGamma_0(D),\chi),$  we have  
$$ \langle I_{2n}(f),\, I_{2n}(f) \rangle$$
$$=2^{-4nk-4n^2-4n+2} D^{2nk+5n^2-3n/2-1/2}\eta_n(f)\displaystyle\prod_{i=1}^{2n} \widetilde\Lambda(i,\, f,\, {\rm Ad},\chi^{i-1}) \prod_{i=2}^{2n} \widetilde{\Lambda}(i,\chi^i),$$
where 
$$\eta_n(f)=\sum_{Q \subset Q_D \atop f_Q=f} \chi_Q((-1)^n).$$
{\rm (2)} Let $m=2n+1$ be a positive odd integer. For a primitive form $f$ in ${\mathfrak S}_{2k}(SL_2({\bf Z})),$  we have  
$$ \langle I_{2n+1}(f),\, I_{2n+1}(f) \rangle$$
$$=2^{-2(2n+1)k-4n^2-6n} D^{2nk+5n^2+5n/2}\displaystyle\prod_{i=1}^{2n+1} \widetilde \Lambda(i,\, f,\, {\rm Ad},\chi^{i-1})\prod_{i=2}^{2n+1}\widetilde{\Lambda}(i,\chi^i).$$
}

\bigskip

\noindent
{\bf Remark.}  In \cite{Ik3} Ikeda showed that $I_m(f)$ is identically zero if and only if $m=2n$ and $\eta_n(f)=0.$ Therefore the above theorem remains valid even if $I_m(f)$ is identically zero.

This type of result was conjectured by Ikeda \cite{Ik3}.  
When $m=2$,  by using the result of Sugano \cite{Su}, Ikeda \cite{Ik3}  has been already proved that
 $$ \langle I_{2}(f),\, I_{2}(f) \rangle=\eta_1(f)2^{-4k-6} D^{2k+3}\widetilde \Lambda(2) \widetilde \Lambda(1,f,{\rm Ad}) \widetilde \Lambda(2,f,{\rm Ad},\chi).$$
His conjecture holds true up to a power of $D.$ In fact, he conjectured that integer powers of $D$  should appear on the right-hand sides of the above formulas. However, half-integer powers of $D$ appear in some cases as shown in the above theorem.

\bigskip

Now put
$${\bf L}(i,f,{\rm Ad},\chi^{i-1})={\widetilde{\Lambda}(i,\, f,\, {\rm Ad},\chi^{i-1}) \over \langle  f,\, f \rangle }$$
 for $i=1,...,m$ 
$${\bf L}(2i,\chi^{2i})=\widetilde{\Lambda}(2i,\chi^{2i}),$$
and
$${\bf L}(2i+1,\chi^{2i+1})=\widetilde{\Lambda}(2i+1,\chi^{2i+1})D^{2i+1/2}$$
for an integer $i \ge 1.$ We note that
$${\bf L}(1,f,{\rm Ad})=
\left\{ \begin{array}{ll}
2^{2k+1} \prod_{q | D}(1+q^{-1}) \ & \ {\rm if} \ f \in {\mathfrak S}_{2k+1}(\varGamma_0(D),\chi) \\
2^{2k}  \ & \ {\rm if} \ f \in {\mathfrak S}_{2k}(SL_2({\bf Z})).
\end{array}
\right.$$
Hence we obtain the following:

\bigskip
\noindent
{\bf Theorem 2.2.} {\it 
Let the notation be as above. Then we have 
 $$\frac{ \langle I_m(f),\, I_m(f) \rangle}{ \langle  f,\, f \rangle ^{m}}=2^{\beta_{n,k}}  \prod_{i=2}^{m} {\bf L}(i,\, f,\, {\rm Ad},\chi^{i-1}) {\bf L}(i,\chi^i)$$
 $$\times \left\{\begin{array}{ll}
 \eta_n(f)D^{2nk+4n^2-n} \prod_{q | D}(1+q^{-1}) & \ {\rm if} \ m=2n \\
 D^{2nk+4n^2+n} & \ {\rm if} \ m=2n+1,
 \end{array}
 \right.$$
 where $\beta_{n,k}$ is an integer depending on $n$ and $k.$
}

\bigskip

It is well known that ${\bf L}(i,\chi^{i})$ is a rational number for any positive integer $i$. Moreover ${\bf L}(i,f,{\rm Ad},\chi^{i-1})$ is an algebraic number and belongs to the Hecke field ${\bf Q}(f)$ for $i=2,....,k'$ where $k'=2k$ or $2k-1$ according as if $m$ is even or odd (cf. Shimura \cite{Sh1}, \cite{Sh2}). Thus we have 

\bigskip

\bigskip

\noindent
{\bf Theorem 2.3.} {\it 
 In addition to the above notation and the assumption, suppose that $m \le 2k$ or $m \le 2k-1$ according as $m$ is even or odd. Then $\displaystyle \frac{ \langle I_m(f),\, I_m(f) \rangle}{ \langle  f,\, f \rangle ^{m}}$  is algebraic, and in particular it belongs to ${\bf Q}(f).$ 
 } 

\bigskip

\section{Rankin-Selberg convolution product }
  To prove Theorem 2.1, we rewrite it in terms of the residue of the Rankin-Selberg convolution product of $I_m(f).$   Let
 $$F(z)=\sum_{A \in \widehat {\mathrm{Her}}_m({\mathcal O})^+} a_F(A){\bf e}({\rm tr}(Az)$$
 be an element of ${\mathfrak S}_{2l}(\varGamma^{(m)},\det^{-l}).$ We then define the Rankin-Selberg series $R(s,F)$ for $F$ by
  $$R(s,F)=\sum_{A \in {\widehat {\mathrm{Her}}_m({\mathcal O})}^+/SL_m(\mathcal O)} {a_{F}(A) \overline{a_{F}(A)} \over (\det A)^{s} e^*(A)},$$
  where  $e^*(A)=\#(\{g \in  SL_m(\mathcal O) \ | \ g^*Ag=A \}).$   
 
  \bigskip
 \noindent
{\bf Proposition 3.1.} {\it
 Put
$$R_m= {2^{2lm+m-1}\prod_{i=2}^m L(i,\chi^{i+1}) \over  D^{m(m-1)/2} \prod_{i=0}^{m-1} L(2m-i,\chi^i)\prod_{i=1}^m \Gamma_{\bf C}(i)\Gamma_{\bf C}(2l-i+1) }.$$ Let $F \in {\mathfrak S}_{2l}(\varGamma^{(m)},\det^{-l}).$  
Then $R(s,F)$ is holomorphic in $s$ for ${\rm Re}(s) > 2l.$ Moreover it can be continued to a meromorphic function on the whole $s$-plane, and has a simple pole at $s=2l$ with the residue $R_m \langle F, F \rangle.$
}

\bigskip
\begin{proof}
The assertion can be proved by a careful analysis of the proof of [\cite{Sh2}, Proposition 22.2]. However, for the convenience of the readers we here give an outline of the  proof. We  define another Rankin-Selberg series $\widetilde {R}(s,F)$ for $F$ by
  $$\widetilde {R}(s,F)=\sum_{A \in {\widehat {\mathrm{Her}}_m({\mathcal O})}^+/GL_m(\mathcal O)} {a_{F}(A) \overline{a_{F}(A)} \over (\det A)^{s}  e(A)},$$
  where  $e(A)=\#(\{g \in  GL_m(\mathcal O) \ | \ g^*Ag=A \}).$ Remark that 
  $$R(s,F)=\#({\mathcal O}^*) \widetilde {R}(s,F).$$
We define the non-holomorphic  Eisenstein 
series 
$E(Z,s)$ for $\varGamma^{(m)}$  by
$$E(Z,s)=(\det Y)^{s}\sum_{M \in \varGamma^{(m)}_{\infty} \backslash \varGamma^{(m)}} |j(M,Z)|^{-2s},$$
where
$\varGamma^{(m)}_{\infty}=\{\mattwo{A}{B}{0}{D} \in \varGamma^{(m)} \}.$  
Then by using the same argument  as in Page 179 of \cite{Sh2}, we obtain
$$\widetilde R(s,F)={1 \over \#({\mathcal O}^*) {\rm vol}({\mathrm{Her}}_m({\bf C})/{\mathrm{Her}}_m({\mathcal O})) \widetilde {\Gamma}_m(s)(4 \pi )^{-ms}}$$
$$ \times \int_{\varGamma^{(m)} \backslash {\mathfrak H}_{m}}  F(Z) \overline{F(Z)}   \overline{E(Z,\bar s-2l+m)} (\det Y)^{2l-2m}dXdY,$$
where ${\rm vol}({\mathrm{Her}}_m({\bf C})/{\mathrm{Her}}_m({\mathcal O}))$ is the volume of  ${\mathrm{Her}}_m({\bf C})/{\mathrm{Her}}_m({\mathcal O})$ with respect to the measure $dX,$ and 
$$\widetilde {\Gamma}_m(s)=\pi^{m(m-1)/2}\prod_{i=0}^{m-1} \Gamma(s-i).$$
By [\cite{Sh1},Theorem 19.7], $ E(Z,s-2l+m)$ is holomorphic in $s$ for ${\rm Re}(s)>2l.$ Moreover it has a meromorphic continuation to the whole $s$-plane, and has a simple pole at $s=2l$ with the residue of the following form:
$$\pi^{m^2} \widetilde {\Gamma}_m(m)^{-1}{2^{m(1-m)-1}\prod_{i=2}^m L(i,\chi^{i+1}) \over {\rm vol}({\mathrm{Her}}_m({\bf C})/{\mathrm{Her}}_m({\mathcal O})) \prod_{i=0}^{m-1} L(2m-i,\chi^i)} .$$
We note that 
$${\rm vol}({\mathrm{Her}}_m({\bf C})/{\mathrm{Her}}_m({\mathcal O}))= 2^{m(1-m)/2}D^{m(m-1)/4}.$$
Thus  we prove the assertion.
\end{proof}

\section{Reduction to local computations}
To prove our main result, we give an explicit formula for $R(s,I_m(f)).$  To do this, we reduce the problem to local computations. 
Let $K_p$ and ${\mathcal O}_p$ be as in Notation. Then $K_p$ is a quadratic extension of ${\bf Q}_p$ or $K_p={\bf Q}_p \oplus {\bf Q}_p.$ 
In the former case let ${\mathcal O}_p$ be the ring of integers in $K_p,$  and $f_p$ the exponent of the conductor of $K_p/{\bf Q}_p.$ If $K_p$ is ramified over ${\bf Q}_p,$ put $e_p=f_p -\delta_{2,p},$ where $\delta_{2,p}$ is Kronecker's delta. 
If  $K_p$ is unramified over ${\bf Q}_p,$ put $e_p=0.$ 
In the latter case, put ${\mathcal O}_p={\bf Z}_p \oplus {\bf Z}_p,$  and $e_p=f_p=0.$  
Moreover put  $\widetilde {\mathrm{Her}}_{m}({\mathcal O}_p)=p^{e_p}\widehat {\mathrm{Her}}_{m}({\mathcal O}_p).$ We note that $\widetilde {\mathrm{Her}}_{m}({\mathcal O}_p)={\mathrm{Her}}_{m}({\mathcal O}_p)$ if $K_p$ is not ramified  over ${\bf Q}_p.$ Let $K$ be an imaginary quadratic extension of ${\bf Q}$ with the discriminant $-D.$ We then put $\widetilde D=\prod_{p | D} p^{e_p},$ and
$\widetilde {{\mathrm{Her}}}_{m}({\mathcal O})=\widetilde D \widehat {\mathrm{Her}}_{m}({\mathcal O}).$ 
Now let $m$ and $l$ be positive integers such that $m \ge l.$ Then for an integer  $a$ and  $A \in \widetilde {{\mathrm{Her}}}_{m}({\mathcal O}_p), B \in \widetilde {{\mathrm{Her}}}_{l}({\mathcal O}_p)$ put 
$${\mathcal A}_a(A,B)=\{X \in
M_{ml}({\mathcal O}_p)/p^aM_{ml}({\mathcal O}_p) \ | \ A[X]-B \in p^a\widetilde {{\mathrm{Her}}}_l({\mathcal O}_p) \},$$
and 
$${\mathcal B}_a(A,B)=\{X \in {\mathcal A}_a(A,B) \ | \ 
  {\rm rank}_{{\mathcal O}_p/p{\mathcal O}_p} X=l \}.$$
Suppose that $A$ and  $B$ are non-degenerate. Then the number $p^{a(-2ml+l^2)}\#{\mathcal A}_a(A,B)$ is independent of $a$ if $a$ is sufficiently large. Hence we define the local density $\alpha_p(A,B)$  representing $B$ by $A$ as
$$\alpha_p(A,B)=\lim_{a \rightarrow
\infty}p^{a(-2ml+l^2)}\#{\mathcal A}_a(A,B).$$
Similarly  we can define the primitive local density $\beta_p(A,B)$ as 
 $$\beta_p(A,B)=\lim_{a \rightarrow
\infty}p^{a(-2ml+l^2)}\#{\mathcal B}_a(A,B)$$
 if $A$ is non-degenerate. We remark that the primitive local density $\beta_p(A,B)$ can be defined even if $B$ is not non-degenerate. 
In particular we write $\alpha_p(A)=\alpha_p(A,A).$ 


Let ${\mathcal U}_1$ be the unitary group defined in Section 1. Namely let 
$${\mathcal U}_1=\{ u \in R_{K/{\bf Q}}(GL_1) \ | \ \overline {u} u =1 \}. $$
For an element $T \in {\mathrm{Her}}_m({\mathcal O}_p),$ let
$$\widetilde {U_{p,T}}=\{ \det X \ | \ X \in {\mathcal U}_T(K_p) \cap GL_m({\mathcal O}_p)) \}.$$ Then $\widetilde {U_{p,T}}$ is a subgroup of $U_{1,p}$ of finite index. We then put  
\\$l_{p,T}=[U_{1,p}:\widetilde {U_{p,T})}].$ 
We also put 
$$u_p=\left\{\begin{array}{ll}
( 1+p^{-1} )^{-1} & \ {\rm if} \ K_p/{\bf Q}_p \ {\rm is \ unramified } \\
 (1-p^{-1})^{-1} & \ {\rm if} \ K_p={\bf Q}_p \oplus {\bf Q}_p \\
 2^{-1}         & \  \ {\rm if} \ K_p/{\bf Q}_p \ {\rm is \ ramified. }
 \end{array} \right. $$
For a subset ${\mathcal T}$ of ${\mathcal O}_p$ put
$${\mathrm{Her}}_m({\mathcal T})={\mathrm{Her}}_m({\mathcal O}_p) \cap M_m({\mathcal T}),$$ and for a subset ${\mathcal S}$ of ${\mathcal O}_p$ put
$${\mathrm{Her}}_m({\mathcal S},{\mathcal T})= \{ A \in {\mathrm{Her}}_m({\mathcal T}) \ | \  \det A \in {\mathcal S} \},$$ 
and $\widetilde {\mathrm{Her}}_m({\mathcal S},{\mathcal T})={\mathrm{Her}}_m({\mathcal S},{\mathcal T}) \cap \widetilde {\mathrm{Her}}_m({\mathcal O}_p).$
In particular if ${\mathcal S}$ consists of a single element $d$ we write ${\mathrm{Her}}_m({\mathcal S},{\mathcal T})$ as ${\mathrm{Her}}_m(d,{\mathcal T}),$ and so on. For  $d \in {\bf Z}_{>0}$ we also define the set ${\mathrm{Her}}_m(d,{\mathcal O})^+$ in a similar way. 
For each $T \in {\widetilde {\mathrm{Her}}}_m({\mathcal O}_p)^{\times}$ put 
$$F_p^{(0)}(T,X)=F_p(p^{-{e}_p}T,X)$$ and 
$$\widetilde F_p^{(0)}(T,X)=\widetilde F_p(p^{-{e}_p}T,X).$$
We remark that 
$$\widetilde F_p^{(0)}(T,X)=X^{-{\rm ord}_p(\det T)}X^{e_pm-f_p[m/2]} F_p^{(0)}(T,p^{-m}X^2).$$
For $d \in {\bf Z}_p^{\times}$ put 
$$\lambda_{m,p}(d,X,Y)=  \sum_{A \in \widetilde {\mathrm{Her}}_m(d,{\mathcal O}_p)/SL_{m}({\mathcal O}_p)} {\widetilde F_p^{(0)}(A,X^{-1})\widetilde F_p^{(0)}(A,Y^{-1}) \over u_pl_{p,A}\alpha_p(A)}.$$
An explicit formula for  $\lambda_{m,p}(p^id_0,X,Y)$ will be given in the next section for $d_0 \in {\bf Z}_p^*$ and
$i \ge 0.$ 
 
 \noindent
 {\bf Theorem 4.1.}  {\it 
  Let $f$ be a primitive form in ${\mathfrak S}_{2k+1}(\varGamma_0(D),\chi)$  or in ${\mathfrak S}_{2k}(SL_2({\bf Z}))$ according as $m=2n$ or $2n+1.$ For such an $f$ and a positive integer $d_0$ put
$$a_m(f;d_0)=\prod_p \lambda_{m,p}(d_0,\alpha_p,\overline {\alpha}_p),$$
where $\alpha_p$ is the Satake $p$-parameter of $f.$ 
Moreover put
$$\mu_{m,k,D}=D^{m(s-2k+l_0)+(2k-l_0)[m/2]-m(m+1)/4-1/2}$$
$$ \times 2^{-c_Dm(s-2k-2n)-m+1} \prod_{i=2}^{m} \Gamma_{\bf C}(i),$$
where $l_0=0$ or $1$ according as $m$ is even or odd, and $c_D=1$ or $0$ according as $2$ divides $D$ or not.  Then for ${\rm Re}(s) >>0,$ we have 
$$R(s,I_m(f))=\mu_{m,k,D} \sum_{d_0=1}^{\infty} a_m(f;d_0) d_0^{-s+2k+2n}.$$

}

\begin{proof}
We note that $R(s,I_m(f))$ can be rewritten as 
$$R(s,I_m(f))= \widetilde {D}^{ms} \sum_{T \in \widetilde {\mathrm{Her}}_m({\mathcal O})^+/SL_m({\mathcal O})} {a_{I_m(f)}({\widetilde D}^{-1}T)\overline{a_{I_m(f)}({\widetilde D}^{-1}T)} \over e^*(T) (\det T)^s }.$$
For $T \in \widetilde {\mathrm{Her}}_m({\mathcal O})^+$ the Fourier coefficient $a_{I_m(f)}({\widetilde D}^{-1}T)$ of $I_m(f)$ is uniquely determined by the genus to which $T$ belongs, and can be expressed as
$$|a_{I_m(f)}({\widetilde D}^{-1}T)|^2=(D^{[m/2]}\widetilde{D}^{-m}\det T)^{2k-l_0}\prod_p \widetilde F_p^{(0)}(T,\alpha_p)\widetilde F_p^{(0)}(T,\overline{\alpha_p}).$$
Thus the assertion follows from  [\cite{Kat3}, Corollary to Proposition 3.2 and Proposition 3.3]. (See also the proof of [\cite{Kat3}, Theorem 3.4].)  
\end{proof}

\bigskip

\section{Formal power series associated with local Siegel series}
  Let $K_p$ be a quadratic extension of ${\bf Q}_p,$ and $\varpi=\varpi_p$ and $ \pi=\pi_p$ be  prime elements of $K_p$ and ${\bf Q}_p$, respectively.   If $K_p$ is unramified over ${\bf Q}_p,$ we take $\varpi=\pi=p.$
 If $K_p$ is ramified over ${\bf Q}_p,$ we take $\pi$ so that $\pi=N_{K_p/{\bf Q}_p}(\varpi).$
Let $K_p={\bf Q}_p \oplus {\bf Q}_p.$ Then put $\varpi=\pi=p.$ 
For $d_0 \in {\bf Z}_p^{\times}$ put 
$$\hat H_{m,p}(d_0,X,Y,t)= \sum_{i=0}^{\infty}\lambda_{m,p}^*(p^id_0,X,Y)t^i,$$
where for $d \in {\bf Z}_p^{\times}$ we define $\lambda_{m,p}^*(p^id_0,X,Y)$ as
$$\lambda_{m,p}^*(d,X,Y)=
\sum_{A \in {\widetilde{\mathrm{Her}}}_m(dN_{K_p/{\bf Q}_p}({\mathcal O}_p^*),{\mathcal O}_p)/GL_{m}({\mathcal O}_p)}  {\widetilde F_p^{(0)}(A,X^{-1})\widetilde F_p^{(0)}(A,Y^{-1}) \over \alpha_p(A)}.$$
We note that 
$$\lambda_{m,p}^*(d,X,Y)=
\sum_{A \in {\widetilde{\mathrm{Her}}}_m(dN_{K_p/{\bf Q}_p}({\mathcal O}_p^*),{\mathcal O}_p)/GL_{m}({\mathcal O}_p)}  {\widetilde F_p^{(0)}(A,X)\widetilde F_p^{(0)}(A,Y) \over \alpha_p(A)}.$$
In  Proposition 5.5.1 we will show that  we have  
$$\lambda_{m,p}^*(d,X,Y)=u_p\lambda_{m,p}(d,X,Y)$$
for $d \in {\bf Z}_p^{\times}$ and therefore   
$$\hat H_{m,p}(d_0,X,Y,t)=u_p \sum_{i=0}^{\infty}\lambda_{m,p}(p^id_0,X,Y)t^i.$$
We also define  $H_{m,p}(d_0,X,Y,t)$ as 
$$H_{m,p}(d_0,X,Y,t)= \sum_{i=0}^{\infty}\lambda_{m,p}^*(\pi^id_0,X,Y)t^i.$$
We note that $H_{m,p}(d_0,X,Y,t)=\hat H_{m,p}(d_0,X,Y,t)$ if $K_p$ is unramified over ${\bf Q}_p$ or $K_p={\bf Q}_p \oplus {\bf Q}_p,$ but it is not necessarily the case if $K_p$ is ramified over ${\bf Q}_p.$
In this section, we give  explicit formulas of $H_{m,p}(d_0,X,Y,t)$ for all prime numbers $p$ (cf. Theorems 5.5.2 and 5.5.3),
and therefore explicit formulas for $\hat H_{m,p}(d_0,X,Y,t)$ (cf. Theorem 5.5.4).

From now on we fix a prime number $p.$ Throughout this section we simply write ${\rm ord}_p$ as ${\rm ord}$ and so on if the prime number $p$ is clear from the context. We also write $\nu_{K_p}$ as $\nu.$ We also simply write $\widetilde {\mathrm{Her}}_{m,p}$ instead of $\widetilde{\mathrm{Her}}_m({\mathcal O}_p),$ and so on.  For a $GL_m({\mathcal O}_p)$-stable subset $\mathcal B$ of ${\mathrm{Her}}_m(K_p)$ we 
 simply write $\sum_{T \in \mathcal B}$ instead of $\sum_{T \in \mathcal B/GL_m({\mathcal O}_p)}$ if there is no fear of confusion.

\subsection{Preliminaries}

\noindent
{ }

\bigskip

 Let $m$ be a positive integer. For a non-negative integer $i \le m$ let 
$${\mathcal D}_{m,i}=GL_m({\mathcal O}_p) \mattwo{1_{m-i}}{0}{0}{\varpi 1_i} GL_m({\mathcal O}_p),$$
 and for $W \in  {\mathcal D}_{m,i},$ put $\varPi_p(W)=(-1)^i p^{i(i-1)a/2},$
 where $a=2$ or $1$ according as $K_p$ is unramified over ${\bf Q}_p$ or not.
Let $K_p={\bf Q}_p \oplus {\bf Q}_p.$ Then for a pair $i=(i_1,i_2)$ of non-negative integers such that $i_1,i_2 \le m,$ let 
$${\mathcal D}_{m,i}=GL_m({\mathcal O}_p) \left(\mattwo{1_{m-i_1}}{0}{0}{p 1_{i_1}},\mattwo{1_{m-i_2}}{0}{0}{p 1_{i_2}} \right) GL_m({\mathcal O}_p),$$
 and for $W \in {\mathcal D}_{m,i}$ put $\varPi_p(W)=(-1)^{i_1+i_2} p^{i_1(i_1-1)/2+i_2(i_2-1)/2}.$ In either case $K_p$ is a quadratic extension of ${\bf Q}_p,$ or  $K_p={\bf Q}_p \oplus {\bf Q}_p,$ we put $\varPi_p(W)=0$ for $W \in M_n({\mathcal O}_p^{\times}) \setminus \bigcup_{i=0}^m {\mathcal D}_{m,i}.$ 

For non-degenerate Hermitian matrices $S$ and $T$ of degree $m,$  we put 
$$\alpha_p(S,T;i)=\lim_{e \longrightarrow \infty}p^{-m^2e} {\mathcal A}_e(S,T;i),$$
where
$${\mathcal A}_e(S,T;i)=\{ \bar X \in M_m({\mathcal O}_p)/p^eM_m({\mathcal O}_p) \in {\mathcal A}_e(S,T) \ | \ X \in {\mathcal D}_{m,i} \}.$$

\bigskip
For two elements $A,A' \in {\mathrm{Her}}_m({\mathcal O}_p)$ we simply write $A \sim_{GL_m({\mathcal O}_p)} A'$ as $A \sim A'$ if there is no fear of confusion. 
For a variables $U$ and $q$ put
$$(U,q)_m=\prod_{i-1}^{m}(1-q^{i-1}U), \qquad  \phi_m(q)=(q,q)_m.$$
We note that $\phi_m(q)=\prod_{i=1}^m (1-q^i).$
Moreover for a prime number $p$ put
 $$\phi_{m,p}(q)=\left\{\begin{array}{ll}
 \phi_m(q^2) \ & \ {\rm if} \ K_p/{\bf Q}_p \ {\rm is \ unramified} \\
 \phi_m(q)^2 \ & \ {\rm if} \ K_p={\bf Q}_p \oplus {\bf Q}_p \\
 \phi_m(q)   \ & \ {\rm if} \  K_p/{\bf Q}_p \ {\rm is \ ramified}
 \end{array}
 \right.$$

\bigskip

\noindent
{\bf Lemma 5.1.1.} {\it
 {\rm (1)} Let $\Omega(S,T)=\{w \in M_m({\mathcal O}_p) \ | \ S[w] \sim T \},$ and $\Omega(S,T;i)= \Omega(S,T) \cap  {\mathcal D}_{m,i}. $ Then we have
$${\alpha_p(S,T) \over \alpha_p(T)}=\#(\Omega(S,T)/GL_m({\mathcal O}_p))p^{-m({\rm ord}(\det T)- {\rm ord}(\det S))},$$
and
$${\alpha_p(S,T;i) \over \alpha_p(T)}=\#(\Omega(S,T;i)/GL_m({\mathcal O}_p))p^{-m({\rm ord}(\det T)- {\rm ord}(\det S))} .$$

{\rm (2)} Let $\widetilde \Omega(S,T)=\{w \in M_m({\mathcal O}_p) \ | \ S \sim T[w^{-1}] \},$ and $\widetilde \Omega(S,T;i)= \widetilde \Omega(S,T) \cap {\mathcal D}_{m,i}.$ Then we have  
$${\alpha_p(S,T) \over \alpha_p(S)}=\#(GL_m({\mathcal O}_p) \backslash \widetilde \Omega(S,T)),$$
and $${\alpha_p(S,T;i) \over \alpha_p(S)}=\#(GL_m({\mathcal O}_p) \backslash \widetilde \Omega(S,T;i)).$$
}

\bigskip

\begin{proof} The assertions for ${\alpha_p(S,T) \over \alpha_p(T)}$ and ${\alpha_p(S,T) \over \alpha_p(S)}$ have been proved in [\cite{Kat3}, Lemma 4.1.3].
The assertions for ${\alpha_p(S,T;i) \over \alpha_p(T)}$ and  ${\alpha_p(S,T;i) \over \alpha_p(S)}$ can also be proved in a similar way.
\end{proof}

\bigskip

We define a reduced matrix.
A non-degenerate square matrix $W=(d_{ij})_{m \times m}$ with entries in ${\bf Z}_p$ is said to be reduced if $d_{ii}=p^{e_i}$ with $e_i$ a non-negative integer, $d_{ij}$ is a non-negative integer such that $d_{ij} \le p^{e_j}-1$ for $i <j,$ and $d_{ij}=0$ for $i >j.$ Let $K_p={\bf Q}_p \oplus {\bf Q}_p.$ Then an element $W=(W_1, W_2)$ of $M_m({\mathcal O}_p)^{\times}$ with $W_1,W_2 \in M_m({\bf Z}_p)^{\times}$ is said to be reduced if $W_1$ and $W_2$  are reduced. 
Let $K_p$ be an unramified  quadratic extension of ${\bf Q}_p,$ and $\theta$ be an element of ${\mathcal O}_p$ such that ${\mathcal O}_p={\bf Z}_p + {\bf Z}_p\theta.$ Then a non-degenerate square matrix $W=(d_{ij})_{m \times m}$ with entries in ${\mathcal O}_p$ is said to be reduced if $d_{ii}=p^{e_i}$ with $e_i$ a non-negative integer,  $d_{ij}=d_{ij}^{(1)}+d_{ij}^{(2)}\theta$ with $d_{ij}^{(1)},d_{ij}^{(2)}$ non-negative integers such that  $d_{ij}^{(1)},d_{ij}^{(2)} \le p^{e_j}-1$ for $i <j,$ and $d_{ij}=0$ for $i >j.$ 
Let $K_p$ be a ramified  quadratic extension of ${\bf Q}_p,$ and $\varpi$ be a prime element of $K_p.$ Then a non-degenerate square matrix $W=(d_{ij})_{m \times m}$ with entries in ${\mathcal O}_p$ is said to be reduced if $d_{ii}=\varpi^{e_i}$ with $e_i$ a non-negative integer,  $d_{ij}=d_{ij}^{(1)}+d_{ij}^{(2)}\varpi$ with  $d_{ij}^{(1)}, d_{ij}^{(2)}$ non-negative integers such that  $d_{ij}^{(1)} \le p^{[(e_j+1)/2]}-1, 0 \le d_{ij}^{(2)} \le p^{[e_j/2]}-1$ for $i <j,$ and $d_{ij}=0$ for $i >j.$ 
 In any case, we can take the set of all reduced matrices as a  complete set of representatives of $GL_m({\mathcal O}_p) \backslash M_m({\mathcal O}_p)^{\times}.$ Let $m$ be an integer.  For $B \in {\widetilde{\mathrm{Her}}}_{m}({\mathcal O}_p)$ put
$$\widetilde \Omega(B)=\{W \in GL_m(K_p) \cap M_m({\mathcal O}_p) \ | \ B[W^{-1}] \in \widetilde{\mathrm{Her}}_{m}({\mathcal O}_p) \}.$$
Moreover put $\widetilde \Omega(B,i)=\widetilde \Omega(B) \cap {\mathcal D}_{m.i}.$ Let $r \le m,$ and $\psi_{r,m}$ be the mapping from $GL_{r}(K_p)$ into $GL_{m}(K_p)$ defined by $\psi_{r,m}(W)=1_{m-r} \bot W.$ 

For a subset ${\mathcal T}$ of ${\mathcal O}_p,$ we put 
$${\mathrm{Her}}_m({\mathcal T})_k=\{ A=(a_{ij})  \in {\mathrm{Her}}_m({\mathcal T}) \ | \  a_{ii} \in \pi^k{\bf Z}_p \}.$$
 
 From now on put
 $${\mathrm{Her}}_{m,*}({\mathcal O}_p)=\left\{\begin{array}{ll}
  {\mathrm{Her}}_{m}({\mathcal O}_p)_1 \ & {\rm if} \ p=2 \ {\rm and} \ f_p=3, \\
 {\mathrm{Her}}_{m}(\varpi {\mathcal O}_p)_1 \ & {\rm if} \ p=2 \ {\rm and} \ f_p=2 \\
 {\mathrm{Her}}_{m}({\mathcal O}_p) \ & \ {\rm  otherwise,} 
 \end{array}
 \right.$$
 where $\varpi$ is a prime element of $K_p.$ 
Moreover put  $i_p=0,$ or $1$ according as $p=2$ and $f_2=2,$ or not.
 Suppose that $K_p/{\bf Q}_p$ is unramified  or $K_p={\bf Q}_p \oplus {\bf Q}_p.$  Then an element $B$ of $\widetilde{\mathrm{Her}}_{m}({\mathcal O}_p)$ can be expressed as 
 $B \sim_{GL_m({\mathcal O}_p)} 1_r \bot pB_2$ with some integer $r$ and $B_2 \in {\mathrm{Her}}_{m-r,*}({\mathcal O}_p).$ 
Suppose that $K_p/{\bf Q}_p$ is ramified. For an even positive integer $r$ define $\Theta_r$ by 
$$\Theta_{r}= \overbrace{\mattwo{0}{\varpi^{i_p}}{\overline {\varpi}^{i_p}}{0} \bot ...\bot \mattwo{0}{\varpi^{i_p}}{\overline{\varpi}^{i_p}}{0}}^{r/2},$$
where $\overline{\varpi}$ is the conjugate of $\varpi$ over ${\bf Q}_p.$
 Then an element $B$ of $\widetilde{\mathrm{Her}}_{m}({\mathcal O}_p)$ is expressed as $B \sim_{GL_m({\mathcal O}_p)} \Theta_r \bot \pi^{i_p}B_2$ with some even integer $r$ and $B_2 \in {\mathrm{Her}}_{m-r,*}({\mathcal O}_p).$
For these results, see Jacobowitz \cite{J}.

\bigskip
\noindent
{\bf Lemma 5.1.2.} {\it 

{\rm (1)} Suppose that  $K_p$ is unramified over  ${\bf Q}_p$ or $K_p={\bf Q}_p \oplus {\bf Q}_p.$ Let $B_1 \in {\mathrm{Her}}_{m-n_0}({\mathcal O}_p).$ Then $\psi_{m-n_0,m}$ induces a bijection from $GL_{m-n_0}({\mathcal O}_p) \backslash \widetilde \Omega(pB_1)$ to $GL_{m}({\mathcal O}_p) \backslash \widetilde \Omega(1_{n_0} \bot pB_1),$
which will be also denoted by $\psi_{m-n_0,m}.$ \\
{\rm (2)} Suppose that  $K_p$ is ramified over ${\bf Q}_p$ and that $n_0$ is even. Let $B_1 \in {\mathrm{Her}}_{m-n_0}({\mathcal O}_p).$ Then $\psi_{m-n_0,m}$ induces a bijection from 
$GL_{m-n_0}({\mathcal O}_p) \backslash \widetilde \Omega(\pi^{i_p}B_1)$ to $GL_{m}({\mathcal O}_p) \backslash \widetilde \Omega(\Theta_{n_0} \bot \pi^{i_p}B_1),$
which will be also denoted by $\psi_{m-n_0,m}.$ Here $i_p$ is the integer defined above. \\
{\rm (3)} The assertions remain valid if we replace $\widetilde \Omega(B)$ with $\widetilde \Omega(B,i).$

}
\begin{proof}
The assertions (1) and (2) are due to [\cite{Kat3}, Lemma 4.1.4].  We prove (3). Assume that $K_p$ is unramified over ${\bf Q}_p$ or $K_p={\bf Q}_p \oplus {\bf Q}_p.$ Clearly $\psi_{m-n_0,m}$ is injective. To prove the surjectivity, take a representative $W$ of an element of $GL_{m}({\mathcal O}_p) \backslash \widetilde \Omega(1_{n_0} \bot B_1).$ Without loss of generality we may assume that $W$ is a reduced matrix
with diagonal elements $p^r \ (0 \le r \le 1)$. Since we have $(1_{n_0} \bot B_1)[W^{-1}]  \in  \widetilde{{\mathrm{Her}}_m}({\mathcal O}_p),$ we have 
$W=\mattwo{1_{n_0}}{0}{0}{W_1}$ with  $W_1 \in \widetilde \Omega(B_1,i).$ This proves the assertion. Similarly the assertion holds in the case $K_p$ is ramified over ${\bf Q}_p.$ 

\end{proof}

  
  
   \subsection{Formal power series of Andrianov type}
  
  \noindent
  
  { }
  
  \bigskip
 
For an element $T \in {\widetilde{\mathrm{Her}}}_{m}({\mathcal O}_p),$  we define a polynomial $\widetilde G_p(T,X,t)$ in $X$ and $t$ by
$$\widetilde G_p(T,X,t)=\sum_{i=0}^{m} \sum_{W \in GL_{m}({\mathcal O}_p) \backslash {\mathcal D}_{m,i}}  \varPi_p(W) t^{\nu(\det W)}\widetilde F_p^{(0)}(T[W^{-1}],X).$$
  We also define a polynomial $G_p(T,X)$ in $X$ by
$$G_p(T,X)=\sum_{i=0}^{m} \sum_{W \in GL_{m}({\mathcal O}_p) \backslash {\mathcal D}_{m,i}} (Xp^{m})^{\nu(\det W)}\varPi_p(W)F_p^{(0)}(T[W^{-1}],X).$$
Moreover for an element $T \in {\widetilde{\mathrm{Her}}}_{m,p}$ we define a polynomial $B_p(T,t)$ in $t$ by  
$$B_p(T,t)={\prod_{i=0}^{m-1}(1-\tau_p^{m+i} p^{m+i}t^2) \over G_p(T,t^2)},$$
where $\tau_p^j=1 $ or $\xi_p$ according as $j$ is even or odd. 
 We note that
$$\widetilde G_p(T,X,1)= X^{-{\rm ord}(\det T)}X^{e_pm-f_p[m/2]}G_p(T,Xp^{-m}).$$
Now we recall several results in [\cite{Kat3}].

\bigskip
 \noindent
 {\bf Lemma 5.2.1.} {\it 
 {\rm [\cite{Kat3}, Corollary to Lemma 4.2.2]} {\rm (1)} Suppose that $K_p$ is unramified over ${\bf Q}_p$ or $K_p={\bf Q}_p \oplus {\bf Q}_p.$ Let $T=1_{m-r} \bot pB_1$  with $B_1 \in {\mathrm{Her}}_r({\mathcal O}_p).$ Then we have
  $$G_p(T,Y)= \prod_{i=0}^{r-1}(1-(\xi_p p)^{m+i}Y).$$
  {\rm (2)} Suppose that $K_p$ is ramified over ${\bf Q}_p.$  Let $T=\Theta_{m-r} \bot \pi^{i_p}B_1$  with $B_1 \in {\mathrm{Her}}_{r,*}({\mathcal O}_p).$ Suppose that $m-r$ is even. Then
   $$G_p(T,Y)=\prod_{i=0}^{[(r-2)/2]}(1-p^{2i+2[(m+1)/2]}Y).$$
   }

  \bigskip
 
\noindent
{\bf Lemma 5.2.2.} {\it
{\rm [\cite{Kat3}, Lemma 4.2.3]} Let $B \in {\widetilde{\mathrm{Her}}}_{m}({\mathcal O}_p).$ Then we have
$$F_p^{(0)}(B,X)= \sum_{W \in  GL_{m}({\mathcal O}_p) \backslash \widetilde \Omega(B)} G_p(B[W^{-1}],X)(p^mX)^{\nu(\det W)}.$$}

 \noindent
{\bf Corollary.} {\it {\rm [\cite{Kat3}, Corollary to Lemma 4.2.3]}
Let $B \in {\widetilde{\mathrm{Her}}}_{m}({\mathcal O}_p).$ Then we have
$$\widetilde F_p^{(0)}(B,X)=X^{e_pm-f_p[m/2]}\sum_{B' \in  {\widetilde{\mathrm{Her}}}_{m}({\mathcal O}_p) / GL_{m}({\mathcal O}_p)  } X^{-{\rm ord}(\det B')}{\alpha_p(B',B) \over \alpha_p(B')}$$$$ \times  G_p(B',p^{-m}X^2)X^{{\rm ord}(\det B)-{\rm ord}(\det B')}.$$
}

By Lemma 5.2.1, we easily obtain:

\bigskip
\noindent 
{\bf Lemma 5.2.3.} {\it  
 {\rm (1)} Suppose that $K_p$ is unramified over ${\bf Q}_p$ or $K_p={\bf Q}_p \oplus {\bf Q}_p.$ Let $T=1_{m-r} \bot pB_1$  with $B_1 \in {\mathrm{Her}}_{r}({\mathcal O}_p).$ Then we have
  $$B_p(T,t)= \prod_{i=r}^{m-1}(1-(\xi_p p)^{m+i}t^2).$$
 
 {\rm (2)} Suppose that $K_p$ is ramified over ${\bf Q}_p.$  Let $T=\Theta_{m-r} \bot p^{i_p}B_1$  with $B_1 \in {\mathrm{Her}}_{r,*}({\mathcal O}_p).$ Then
   $$B_p(T,t)=\prod_{i=[(r-1)/2]+1}^{[(m-2)/2]}(1-p^{2i+2[(m+1)/2]}t^2).$$
   }

 \bigskip
  
  For a non-degenerate semi-integral matrix $T$ over ${\mathcal O}_p$ of degree $n,$  put
$$S_p(T,X,t)=\sum_{W \in M_m({\mathcal O}_p)^{\times}/GL_m({\mathcal O}_p)} \widetilde F_p^{(0)}(T[W],X)t^{\nu(\det W)}.$$
This type of formal power series  was first introduced by Andrianov [A] to study the standard $L$-functions of Siegel modular forms of integral weight. Thus we call it the formal power series of Andrianov type. (See also \cite{Bo}, \cite{K-K2}).  The following proposition can easily be proved by (1) of Lemma 5.1.1.

  \bigskip
  \noindent
{\bf Proposition 5.2.4.} {\it 
  Let $T \in {\widetilde{\mathrm{Her}}}_{m}({\mathcal O}_p).$ Then we have  $$\sum_{B \in \widetilde{\mathrm{Her}}_{m}({\mathcal O}_p)}{\widetilde F_p^{(0)}(B,X)\alpha_p(T,B) \over \alpha_p(B)}t^{{\rm ord}(\det B)}= t^{{\rm ord}(\det T)}S_p(T,X,p^{-m}t).$$
  }

  \bigskip
 
 Put ${\mathcal K}^{(m)}={\mathcal K}_0^{(m)}{\mathcal U}^{(m)}({\bf R}).$ Let ${\mathcal H}({\mathcal U}^{(m)}({\bf A}),{\mathcal K}^{(m)})$ be the Hecke ring associated with the Hecke pair $({\mathcal U}^{(m)}({\bf A}),{\mathcal K}^{(m)}).$ 
Then  ${\mathcal H}({\mathcal U}^{(m)}({\bf A}),{\mathcal K}^{(m)})$ acts on \\
 ${\mathcal M}_{2l}({\mathcal U}^{(m)}({\bf Q}) \backslash  {\mathcal U}^{(m)}({\bf A}),\det^{-l})$ as in \cite {Ik3}. We call an element $F$ of \\
${\mathcal M}_{2l}({\mathcal U}^{(m)}({\bf Q}) \backslash  {\mathcal U}^{(m)}({\bf A}),\det^{-l})$ a Hecke eigenform if it is a common eigenfunction of all Hecke operators $T$ in  ${\mathcal H}({\mathcal U}^{(m)}({\bf A}),{\mathcal K}^{(m)}).$ Then for each element \\ $r \in GL_m({\bf A}) \cap \prod_p M_m({\mathcal O}_p),$ let $\lambda_F(r)$ be the eigenvalue of 
${\mathcal K}^{(m)}\mattwo{r^{-1}}{0}{0}{r^*}{\mathcal K}^{(m)}$ with respect to $F,$ and
 define a Dirichlet series ${\mathfrak T}(s,F)$ by
  $${\mathfrak T}(s,F)=\sum_{r \in {\mathcal K}^{(m)}\backslash (GL_m({\bf A}) \cap \prod_p M_m({\mathcal O}_p))/{\mathcal K}^{(m)}} \lambda_F(r) |\det r|_{\bf A}^s,$$
  where 
  $|\det r|_{\bf A}=\prod_p |\det r_p|_{K_p}$ for $r=(r_p) \in GL_m({\bf A}) \cap \prod_p M_m({\mathcal O}_p).$ Then there exists an Euler product ${\mathcal Z}(s,F)$ such that
  $${\mathfrak T}(s,F)=\prod_{i=1}^m L(2s-i+1,\chi^{i-1}){\mathcal Z}(s,F).$$
 We then put 
 $${\mathcal L}(s,F,{\rm st})={\mathcal Z}(s+m-1/2,F),$$ and call it the standard $L$-function of $F$ in the sense of Shimura. We note that our standard $L$-function coincides with that in \cite{Ik3} up to Euler factors at ramified primes. 

Now we define the Eisenstein series on  ${\mathcal U}^{(m)}({\bf A})$ and consider its standard $L$-function in the sense of Shimura. 
Let ${\mathcal P}$ be the maximal parabolic subgroup of ${\mathcal U}^{(m,m)}$ defined by
$${\mathcal P}(R)=\{\gamma =\smallmattwo{a}{b}{0}{d} \in {\mathcal U}^{(m,m)}(R)\}$$
 for any ${\bf Q}$-algebra $R.$ 
Write an element $g=(g_v) \in {\mathcal U}^{(m)}({\bf A})$ as
$$(g_p)_{p <\infty}=\left(\smallmattwo{a_p}{b_p}{0}{d_p}\right)_{p < \infty} (\kappa_p)_{p <\infty}$$
with $\left(\smallmattwo{a_p}{b_p}{0}{d_p}\right)_{p < \infty} \in \prod_{p <\infty} {\mathcal P}({\bf Q}_p)$ and $(\kappa_p)_{p <\infty}  \in {\mathcal K}_0,$ 
and  define the function on ${\mathcal U^{(m)}}({\bf A})$   by
$${\bf f}_{2l}(g)=\prod_p|\det (d_p\overline {d}_p)|_p^{-l}j(g_{\infty},{\bf i})^{-2l} (\det g_{\infty})^l.$$
Let $l$ be a integer such that $l > m$. We then define the normalized Eisenstein series as 
$${\bf E}_{2l}^{(m)}(g)=2^{-m}\prod_{i=1}^m L(i-2l,\chi^{i-1})\sum_{\gamma \in {\mathcal P}({\bf Q}) \backslash {\mathcal U}^{(m)}({\bf Q})} {\bf f}_{2l}(\gamma g).$$
Put
$${\mathcal E}_{2l,m}^{(i)}(Z)=2^{-m}\prod_{j=1}^{m} L(j-2l,\chi^{j-1}) $$
$$ \times \prod_p|\det (t_{i,p})\det (\overline{t_{i,p}})|_p^{l} \sum_{g \in (\varGamma_i \cap {\mathcal P}({\bf Q})) \backslash \varGamma_i} (\det g)^l j(g,Z)^{-2l}$$
for $i=1,\ldots,h,$ where $(t_{i,p})$ be the element of ${\bf G}^{(m)}({\bf A}_f)$ defined in Section 2. Then ${\bf E}_{2l}^{(m)}$ is written as
$${\bf E}_{2l}^{(m)}=({\mathcal E}_{2l,m}^{(1)},{\mathcal E}_{2l,m}^{(2)},\ldots,{\mathcal E}_{2l,m}^{(h)})^{\sharp}.$$

Now put
$${\mathcal L}_{m,p}(X,t)$$
$$= \left\{ 
\begin{array}{ll}
\displaystyle \prod_{i=1}^m\{(1-p^{-m+2i-1}X^2t^2)(1-p^{-m+2i-1}X^{-2}t^2)\}^{-1} & {\rm if} \ K_p/{\bf Q}_p \ {\rm  is \ unramified }\\
\displaystyle  \prod_{i=1}^m\{(1-p^{-m/2+i-1/2}Xt)^2(1-p^{-m/2+i-1/2}X^{-1}t)^2\}^{-1} & {\rm if} \ K_p={\bf Q}_p \oplus {\bf Q}_p \\
\displaystyle  \prod_{i=1}^m \{(1-p^{-m/2+i-1/2}Xt)(1-p^{-m/2+i-1/2}X^{-1}t)\}^{-1} & {\rm if}  \ K_p/{\bf Q}_p \ {\rm  is \  ramified.} 
\end{array} 
\right.
 $$
 
 \bigskip
 
\noindent
 {\bf Proposition 5.2.5.} {\it
 ${\bf E}_{2l}^{(m)}$ is a Hecke eigenform in ${\mathcal M}_{2l}({\mathcal U}^{(m)}({\bf Q}) \backslash  {\mathcal U}^{(m)}({\bf A}),\det^{-l}),$ and its standard $L$-function ${\mathcal L}(s,{\bf E}_{2l}^{(m)},{\rm st})$ in the sense of Shimura is given by 
$${\mathcal L}(s,{\bf E}_{2l}^{(m)},{\rm st})=\prod_p {\mathcal L}_{m,p}(p^{-l+m/2},p^{-s}).$$
  }
\begin{proof} The assertion is more or less well known (cf. [\cite{Ik3}, Proposition 13.5]). But for the sake of completeness, we here give an outline of the proof. For each prime number $p$ let ${\mathcal K}_p^{(m)}={\mathcal U}_m({\bf Q}_p) \cap GL_{2m}({\mathcal O}_p)$. Moreover, for each $\eta \in  {\mathcal U}_m({\bf Q}_p)$ we write $\eta = \left(\begin{matrix} a_{\eta} & b_{\eta} \\ c_{\eta} & d_{\eta} \end{matrix}\right)$ with 
$a_{\eta}, b_{\eta},c_{\eta}$ and $d_{\eta} \in M_m(K_p)$.  First assume that $K_p$ is a field. Then for any $u \in {\mathcal U}_m({\bf Q}_p)$, we can write the coset ${\mathcal K}_p^{(m)}u{\mathcal K}_p^{(m)}$ as 
$${\mathcal K}_p^{(m)} u{\mathcal K}_p^{(m)}=\bigsqcup_{\eta} {\mathcal K}_p^{(m)} \left(\begin{matrix} a_{\eta} & b_{\eta} \\ 0 & d_{\eta} \end{matrix}\right),$$
where $d_{\eta}$ is an upper triangular matrix whose diagonal components are $\varpi^{e_1(\eta)},\ldots, \varpi^{e_m(\eta)}$ with $e_1(\eta),\ldots,e_m(\eta) \in {\bf Z}$. Then, by a simple computation we have
$${\bf E}_{2l}^{(m)}|{\mathcal K}_p^{(m)}u{\mathcal K}_p^{(m)} =\sum_{\eta} q^{-l(e_1(\eta)+\cdots+e_m(\eta))} {\bf E}_{2l}^{(m)},$$
where $q=p^2$ or $p$ according as $K_p/{\bf Q}_p$ is unramified or ramified. 
 We note that $q^{-l(e_1(\eta)+\cdots+e_m(\eta))}=\prod_{i=1}^m (q^{-i} q^{-l+i})^{e_i(\eta)}$. Thus, by
[\cite{Sh1}, (16.1.3)], [\cite{Sh2}, Theorem 19.8] and [\cite{Sh2}, 20.6], we can prove that the Euler factor of ${\mathcal L}(s,{\bf E}_{2l}^{(m)},{\rm st})$ at $p$ is 
${\mathcal L}_{m,p}(p^{-l+m/2},p^{-s})$.  Next assume that $K_p={\bf Q}_p \oplus {\bf Q}_p$. Then, by [\cite{Sh2}, p. 163], 
for any $u \in {\mathcal U}_m({\bf Q}_p)$, we can write the coset ${\mathcal K}_p^{(m)}u{\mathcal K}_p^{(m)}$ as 
$${\mathcal K}_p^{(m)} u{\mathcal K}_p^{(m)}=\bigsqcup_{\eta} {\mathcal K}_p^{(m)} \left(\begin{matrix} a_{\eta} & b_{\eta} \\ 0 & d_{\eta} \end{matrix}\right),$$
where $d_{\eta}$ is a pair of upper triangular matrices  whose diagonal components are $p^{e_1(\eta)},\ldots, p^{e_m(\eta)}$ with $e_1(\eta),\ldots,e_m(\eta) \in {\bf Z}$ and  $p^{e_{m+1}(\eta)},\ldots, p^{e_{2m}(\eta)}$ with $e_{m+1}(\eta),\ldots,e_{2m}(\eta) \in {\bf Z}$, respectively. Then, by a simple computation we have
$${\bf E}_{2l}^{(m)}|{\mathcal K}_p^{(m)}u{\mathcal K}_p^{(m)} =\sum_{\eta} p^{-l(e_1(\eta)+\cdots+e_{2m}(\eta))} {\bf E}_{2l}^{(m)}.$$
We note that $p^{-l(e_1(\eta)+\cdots+e_{2m}(\eta))}=\prod_{i=1}^m (p^{-i} p^{-l+i})^{e_i(\eta)}(p^{-i}p^{-l+i})^{e_{m+i}(\eta)}$. Thus, by
[\cite{Sh2}, p. 163], [\cite{Sh2}, Theorem 19.8] and [\cite{Sh2}, 20.6], we can also prove that the Euler factor of ${\mathcal L}(s,{\bf E}_{2l}^{(m)},{\rm st})$ at $p$ is 
${\mathcal L}_{m,p}(p^{-l+m/2},p^{-s})$. This completes the proof.

 
\end{proof}
 
 \bigskip
For an element $x =(x_v) \in {\bf A}$ put 
${\bf e}_{\bf A}(x)={\bf e}(x_{\infty})\prod_{p<\infty} {\bf e}_p(-x_p)$. We also denote by ${\mathcal{HER}}_m$ the algebraic group defined over ${\bf Q}$ such that ${\mathcal{HER}}_m(S)={\mathrm{Her}}_m(S \otimes K)$ for any ${\bf Q}$-algebra $S$.  Then 
for any $u \in G_m({\bf A})$ and $s \in {\mathcal{HER}}_m({\bf A})$ we have the following Fourier expansion:
$${\bf E}_{2l}^{(m)}\left(\left(\begin{smallmatrix} u & (u^*)^{-1} s \\ 0 & (u^*)^{-1}\end{smallmatrix}\right)\right)=
(\det u  \ \overline{\det u})^l \sum_{T \in {\mathrm{Her}}_m(K)} c_{2l}^{(m)}(T;u) {\bf e}(\sqrt{-1} {\rm tr}(u^*Tu)){\bf e}_{\bf A}({\rm tr}(As)),$$
where $c_{2l}^{(m)}(T;u)$ is a complex number depending only on ${\bf E}_{2l}^{(m)}, T , (u_p)_{p <\infty}$ and $(uu^*)_{\infty}$ (cf. [\cite{Sh1}, Proposition 18.3). Here we have $c_{2l}^{(m)}(T;u) \not=0$ only if $T$ is semi-positive definite.

\noindent
{\bf Remark.} For any $T \in  {\mathrm{Her}}_m(K)^+,$ the $T$-th Fourier coefficient $c_{2l,m}^{(i)}(T)$ of  ${\mathcal E}_{2l,m}^{(i)}(Z)$ is equal to $c_{2l}^{(m)}(T,(t_{i,p}))$ (cf. [\cite{Sh2}, (20.9f)]), and it is  given by 
$$ A_m|\gamma (T)|^{l-m/2} \prod_p |\det (t_{i,p}) |_{K_p}^{m/2}\widetilde F_p(t_{i,p}^*T t_{i,p},p^{-l+m/2}),$$
where $A_m=(-1)^m$ or $1$ according as $m=2n$ or $m=2n+1$
 (cf. \cite{Ik2}, pages 1134-1135).  We notice that $A_m$ appears in the above formula because the definition of $\widetilde F_p(*,X)$ is a slightly different from that in \cite{Ik2} as remarked in Section 2. In general, for any $T \in  {\mathrm{Her}}_m(K)^+$ and $u=(u_p) \in {\bf G}^{(m)}({\bf A}_f)$ we have
 $$c_{2l}^{(m)}(T;u)=A_m|\gamma (T)|^{l-m/2} \prod_p |\det u_p |_{K_p}^{m/2}\widetilde F_p(u_p^*T u_p,p^{-l+m/2}).$$ 
This can be proved in the same way as above.

\noindent
{\bf Theorem 5.2.6.} {\it 
 Let $T$ be an element of $\widetilde {\mathrm{Her}}_{m}({\mathcal O}_p)^{\times}.$ Then we have
$$S_p(T,X,t)=B_p(T,p^{-m/2}t)\widetilde  G_p(T,X,t) {\mathcal L}_{m,p}(X,p^{m/2-1/2}t).$$
}

\begin{proof}
Take an element $\widetilde T \in \widetilde {\mathrm{Her}}_m({\mathcal O})^+$ such that $\widetilde T \sim_{GL_m({\mathcal O}_p)} T.$ 
Then we have 
$$S_p(\widetilde T,X,t)=S_p(T,X,t)$$ 
and 
$$B_p(\widetilde T,p^{-m/2}t)\widetilde  G_p(\widetilde T,X,t)=B_p(T,p^{-m/2}t)\widetilde  G_p(T,X,t).$$
Write $S_p(\widetilde T,X,t)$ and $B_p(\widetilde T,p^{-m/2}t)\widetilde G_p(\widetilde T,X,t) {\mathcal L}_{m,p}(X,p^{m/2-1/2}t)$ as
$$S_p(\widetilde T,X,t)=\sum_{i=0}^{\infty} r_i(X)t^i,$$
and
$$B_p(\widetilde T,p^{-m/2}t)\widetilde G_p(\widetilde T,X,t) {\mathcal L}_{m,p}(X,p^{m/2-1/2}t)=\sum_{i=0}^{\infty} s_i(X)t^i.$$
Then $r_i(X)$ and $s_i(X)$ are  polynomials in $X$ and $X^{-1}.$ 
 For a positive integer $l$ and $A \in \widehat {\mathrm{Her}}_m({\mathcal O})^+,$ 
put
$$D_p(s,A,{\bf E}_{2l}^{(m)})=\sum_{W \in M_m({\mathcal O}_p)^{\times}/GL_m({\mathcal O}_p)} |\det W|_{K_p}^{-m}c_{2l}^{(m)}(A,\widetilde W) p^{-s\nu_{K_p}(\det W)},$$
and 
$$\widetilde G_{2l,m}(A,s) =\sum_{W \in GL_m({\mathcal O}_p) \backslash M_m({\mathcal O}_p)^{\times}}  \varPi_p(W) c_{2l}^{(m)}(A, \widetilde W^{-1})p^{-s \nu_{K_p}(\det W)},$$
where for $V \in  M_m({\mathcal K}_p)^{\times}$ we denote by $\widetilde V=(V_q)$  the element of ${\bf G}^{(m)}({\bf A}_f)$ such that $V_p=V$ and $V_q=1_m$ for any $q \not=p$.  
 Then by Proposition 5.2.5 and by using the same argument as in the proof of  [\cite{Sh2}, Theorem 20.7], we obtain
$$D_p(s+m/2,\widetilde D^{-1} \widetilde T,{\bf E}_{2l}^{(m)})$$
$$=\widetilde G_{2l,m}(\widetilde D^{-1} \widetilde T,s+m/2) B_p(\widetilde T,p^{-s-m/2}) {\mathcal L}_{m,p}(p^{-l+m/2},p^{m/2-1/2-s})$$
for any positive integer $l >m$. By the above remark, for any $A \in {\mathrm{Her}}_m(K)^+$ and  $V \in M_m({\mathcal K}_p)^{\times}$ we have
$$c_{2l}^{(m)}(A,\widetilde V)=d(l,m;A) |\det V|_{K_p}^{m/2} \widetilde F_p(V^*AV,p^{-l+m/2}) ,$$
where $d(l,m;A)=A_m|\gamma(A)|^{l-m/2}\prod_{q \not=p} \widetilde F_q(A,q^{-l+m/2})$.
Hence we have 
$$D_p(s+m/2,\widetilde D^{-1} \widetilde T,{\bf E}_{2l}^{(m)})=d(l,m;\widetilde D^{-1}\widetilde T) S_p(\widetilde T,p^{-l+m/2},p^{-s}),$$
and 
$$\widetilde G_{2l,m}(\widetilde D^{-1} \widetilde T,s+m/2)=d(l,m;\widetilde D^{-1}\widetilde T)\widetilde G_p(\widetilde T,p^{-l+m/2},p^{-s}) ,$$
and therefore
$$d(l,m;\widetilde D^{-1} \widetilde T) S_p(\widetilde T,p^{-l+m/2},p^{-s})$$
$$=d(l,m;\widetilde D^{-1} \widetilde T) B_p(\widetilde T,p^{-s-m/2})\widetilde G_p(\widetilde T,p^{-l+m/2},p^{-s}) {\mathcal L}_{m,p}(p^{-l+m/2},p^{m/2-1/2-s})$$
for any positive integer $l >m.$  We note that $d(l,m;\widetilde D^{-1} \widetilde T) \not=0$ for $l >m$. Hence we have
$$  S_p(\widetilde T,p^{-l+m/2},t)=B_p(\widetilde T,p^{-m/2}t)\widetilde G_p(\widetilde T,p^{-l+m/2},t) {\mathcal L}_{m,p}(p^{-l+m/2},p^{m/2-1/2}t)$$
for any integer $l >m$. This implies that $r_i(p^{-l+m/2})=s_i(p^{-l+m/2})$ for infinitely many positive integers $l.$  Hence we have $r_i(X)=s_i(X).$
\end{proof}

\bigskip
 
  Now by Theorem 5.2.6, we can rewrite $H_{m,p}(d_0,X,Y,t)$ in terms of $G_p(B',Y),  B_p(T,t)$ and $\widetilde G_p(T,X,t)$ in the following way:  For $d_0 \in {\bf Z}_p^{\times}$  put
  $$\widetilde {\mathcal F}_{m,p}(d_0)=\bigcup_{i=0}^{\infty} \widetilde {\mathrm{Her}}_m(\pi^id_0N_{K_p/{\bf Q}_p}({\mathcal O}_p^*),{\mathcal O}_p),$$
  and   define a  formal power series $R_m(d_0,X,Y,t)$ in $t$ by 
\begin{align*} R_m(d_0,X,Y,t)=\sum_{B' \in \widetilde {\mathcal F}_{m,p}(d_0)} { \widetilde G_p(B',X,p^{-m}Yt)   \over \alpha_p(B') } \\
 \times (tY^{-1})^{{\rm ord}(\det B')} B_p(B',p^{-3m/2}Yt) G_p(B',p^{-m}Y^2).
\end{align*}
 
 \bigskip
 \noindent
 {\bf Theorem 5.2.7.} {\it
  We have 
 $$H_{m,p}(d_0,X,Y,t)=Y^{e_pm-f_p[m/2]}R_{m,p}(d_0,X,Y,t) {\mathcal L}_{m,p}(X,tYp^{-m/2-1/2})$$
  for $d_0 \in {\bf Z}_p^{\times}.$ 
 } 

 \begin{proof}  We note that $H_{m,p}(d_0,X,Y,t)$ can be written as
 $$H_{m,p}(d_0,X,Y,t)=\sum_{B \in \widetilde {\mathcal F}_{m,p}(d_0)} t^{{\rm ord}(\det B)} {\widetilde F_p^{(0)}(B,X)\widetilde F_p^{(0)}(B,Y) \over \alpha_p(B)}.$$
 Hence by Corollary to Lemma 5.2.2, we have
 $$H_{m,p}(d_0,X,Y,t)=Y^{e_pm-f_p[m/2]}\sum_{B \in \widetilde {\mathcal F}_{m,p}(d_0)}{t^{{\rm ord}(\det B)} \widetilde F_p^{(0)}(B,X) \over \alpha_p(B)} $$
$$ \times \sum_{B' \in \widetilde{\mathrm{Her}}_{m}({\mathcal O}_p)} {Y^{-{\rm ord}(\det {B'})}G_p(B',p^{-m}Y^2) \alpha_p(B',B) \over \alpha_p(B')} Y^{{\rm ord}(\det B)-{\rm ord}(\det B')}.$$
Let $B,B' \in \widetilde{\mathrm{Her}}_{m}({\mathcal O}_p),$ and suppose that $\alpha_p(B',B) \not=0.$ Then we note that $B \in \widetilde {\mathcal F}_{m,p}(d_0)$ if and only if $B' \in \widetilde{\mathcal F}_{m,p}(d_0).$ Hence by Proposition 5.2.4 and Theorem 5.2.6 we have 
$$Y^{-e_pm+f_p[m/2]}H_{m,p}(d_0,X,Y,t)= \sum_{B' \in {\widetilde{\mathcal F}}_{m,p}(d_0)} {G_p(B',p^{-m}Y^2)Y^{-2{\rm ord}(\det B')} \over \alpha_p(B') } $$
$$ \times \sum_{B \in \widetilde {\mathrm{Her}}_{m}({\mathcal O}_p)} {\widetilde F_p^{(0)}(B,X) \alpha_p(B',B) \over \alpha_p(B)}(tY)^{{\rm ord}(\det B)}$$
$$=\sum_{B' \in {\widetilde{\mathcal F}}_{m,p}(d_0)} { G_p(B',p^{-m}Y^2)Y^{-2{\rm ord}(\det B')} \over \alpha_p(B') } (tY)^{{\rm ord}(\det B')}S_p(B',X,tYp^{-m}) $$
$$=\sum_{B' \in {\widetilde{\mathcal F}}_{m,p}(d_0)} { \widetilde G_p(B',X,p^{-m}Yt)   \over \alpha_p(B') } (tY^{-1})^{{\rm ord}(\det B')}$$
$$\times B_p(B',p^{-3m/2}Yt) G_p(B',p^{-m}Y^2) {\mathcal L}_{m,p}(X,tYp^{-m/2-1/2}).$$ 
\end{proof}
 
    \subsection{Formal power series of modified Koecher-Maass type} 
  
  \noindent
  { }
  
  \bigskip

  Let $r$ be a positive  integer, and  $d_{0} \in {\bf Z}_p^*.$ We then  define a  formal power series $P_{r}(d_0,X,t)$ in $t$ by 
$$P_{r}(d_0,X,t)=\sum_{B \in \widetilde{\mathcal F}_{r,p}(d_0)} {\widetilde F_p^{(0)}(B,X) \over \alpha_p(B)}t^{{\rm ord}(\det B)}.$$
 This type of formal power series appears in an explicit formula of the Koecher-Maass series associated with 
 the Siegel Eisenstein series and the Ikeda lift (cf. \cite{I-K2}, \cite{I-K3}). Thus we call this the formal power series of Koecher-Maass type. To prove Theorems 5.5.1 and 5.5.2, the main results of Section 5, we define a formal power series $\widetilde P_{r}(d_0,X,Y,t)$ in $t$ by 
 $$\widetilde P_{r}(d_0,X,Y,t)= \sum_{B' \in \widetilde{\mathcal F}_{r,p}(d_0)} {\widetilde G_p(B',X,t Y) \over \alpha_p(B') }  (tY^{-1})^{{\rm ord}(\det B')}.$$
   The relation between $\widetilde P_{r}(d_0,X,Y,t)$ and
 $P_{r}(d_0,X,t)$ will be given in the following proposition:

  \bigskip

  \bigskip
\noindent
{\bf Proposition 5.3.1.} {\it
  
\noindent
{\rm (1)} Suppose that $K_p$ is unramified over ${\bf Q}_p.$ 
Then 
$$\widetilde P_{r}(d_0,X,Y,t)= P_{r}(d_0,X,tY^{-1}) \prod_{i=1}^{r} (1-t^4p^{-2r-2+2i}).$$

\noindent
{\rm (2)} Suppose that $K_p={\bf Q}_p \oplus {\bf Q}_p.$ 
Then 
$$\widetilde P_{r}(d_0,X,Y,t)= P_{r}(d_0,X,tY^{-1}) \prod_{i=1}^{r} (1-t^2p^{-r-1+i})^2.$$

\noindent
{\rm (3)} Suppose that $K_p$ is ramified over ${\bf Q}_p.$ Then 
$$\widetilde P_{r}(d_0,X,Y,t)= P_{r}(d_0,X,tY^{-1}) \prod_{i=1}^{r} (1-t^2p^{-r-1+i}).$$
}

\begin{proof} First suppose that  $K_p$ is a quadratic extension of ${\bf Q}_p.$ For each non-negative integer $i \le r$ put 
$$P_{r,i}(d_0,X,t)=\sum_{B \in {\widetilde{\mathcal F}}_{r,p}(d_0)} \sum_{W \in GL_r({\mathcal O}_p) \backslash {\mathcal D}_{r,i}} {\widetilde F_p^{(0)}(B[W^{-1}],X) \over \alpha_p(B)}t^{{\rm ord}(\det B)}.$$
Then by (2) of Lemma 5.1.1 we have
$$P_{r,i}(d_0,X,t)=\sum_{B \in {\widetilde{\mathcal F}}_{r,p}(d_0)} {1 \over \alpha_p(B)} \sum_{B' \in \widetilde{\mathrm{Her}}_{r}({\mathcal O}_p)} {\widetilde F_p^{(0)}(B',X) \alpha_p(B',B;i) \over  \alpha_p(B')} t^{{\rm ord}(\det B)}.$$
Let $B,B' \in \widetilde{\mathrm{Her}}_{r}({\mathcal O}_p),$ and suppose that $\alpha_p(B',B;i) \not=0.$ Then we note that $B \in \widetilde{\mathcal F}_{r,p}(d_0)$ if and only if $B' \in \widetilde{\mathcal F}_{r,p}(d_0).$ Thus by (1) of Lemma 5.1.1 we have
\begin{align*}
& P_{r,i}(d_0,X,t)\\
&=\sum_{B' \in \widetilde{\mathcal F}_{r,p}(d_0) } {\widetilde F_p^{(0)}(B',X)  \over  \alpha_p(B')}  \sum_{B \in \widetilde{\mathrm{Her}}_{r}({\mathcal O}_p)}  t^{{\rm ord}(\det B)} {\alpha_p(B',B;i) \over \alpha_p(B)}\\
&=\sum_{B' \in \widetilde{\mathcal F}_{r,p}(d_0) } {\widetilde F_p^{(0)}(B',X)  \over  \alpha_p(B')} t^{{\rm ord}(\det B')} \# ({\mathcal D}_{r,i}/GL_r({\mathcal O}_p))(tp^{-r})^{ei},
\end{align*}
where $e=2$ or $1$ according as $K_p/{\bf Q}_p$ is unramified or ramified.
By using the same argument as in the proof of Lemma 3.2.18 of Andrianov \cite {A}, we have
$$\#({\mathcal D}_{r,i}/GL_r({\mathcal O}_p))={\phi_{r}(p^e) \over \phi_i(p^e)\phi_{r-i}(p^e)}.$$
Hence we have
\begin{align*}
&P_{r,i}(d_0,X,t) \\
&=\sum_{B'  \in \widetilde{\mathcal F}_{r,p}(d_0) } {\widetilde F_p^{(0)}(B',X)  \over  \alpha_p(B')} t^{{\rm ord}(\det B')}
{\phi_{r}(p^e) \over \phi_i(p^e)\phi_{r-i}(p^e)} (tp^{-r})^{ei}\\
&={\phi_{r}(p^e) \over \phi_i(p^e)\phi_{r-i}(p^e)}  P_{r}(d_0,X,t)(tp^{-r})^{ei}.
\end{align*}
 Then  we have 
$$ \widetilde P_{r}(d_0,X,Y,t)=\sum_{i=0}^{r} (-1)^i p^{i(i-1)e/2}(tY)^{ei}  P_{r,i}(d_0,X,tY^{-1}) .$$
Hence we have 
\begin{align*}
& \widetilde P_{r}(d_0,X,Y,t)=\sum_{i=0}^{r}(-1)^i p^{i(i+1)e/2}(p^{e(-r-1)}t^{2e})^{i} {\phi_{r}(p^e) \over \phi_i(p^e) \phi_{r-i}(p^e)} P_{r}(d_0,X,tY^{-1})\\
&=P_{r}(d_0,X,tY^{-1}) \prod_{i=1}^{r} (1-t^{2e}p^{e(-r-1+i)}).
\end{align*}
Next suppose that $K_p={\bf Q}_p \oplus {\bf Q}_p.$ For a pair $i=(i_1,i_2)$ of non-negative integers such that $i_1,i_2 \le r,$ put 
$$P_{r,i}(d_0,X,t)=\sum_{B \in \widetilde{\mathcal F}_{r,p}(d_0)} \sum_{W \in  GL_r({\mathcal O}_p) \backslash  {\mathcal D}_{r,i}} {\widetilde F_p^{(0)}(B[W^{-1}],X) \over \alpha_p(B)}t^{{\rm ord}(\det B)}.$$
Then by using the same argument as above we can prove that
$$P_{r,i}(d_0,X,t)={\phi_{r}(p) \over \phi_{i_1}(p)\phi_{r-i_1}(p)} {\phi_{r}(p) \over \phi_{i_2}(p)\phi_{r-i_2}(p)}P_{r}(d_0,X,t)(tp^{-r})^{i_1+i_2}.$$
Hence we have 
\begin{align*}
& \widetilde P_{r}(d_0,X,Y,t)\\
&=\sum_{i_1=0}^{r}\sum_{i_2=0}^{r} (-1)^{i_1+i_2} p^{i_1(i_1+1)/2+i_2(i_2+1)/2}(p^{-r-1}t^2)^{i_1+i_2} \\
& \times {\phi_{r}(p) \over \phi_{i_1}(p)\phi_{r-i_1}(p)} {\phi_{r}(p) \over \phi_{i_2}(p)\phi_{r-i_2}(p)} P_{r}(d_0,X,tY^{-1})\\
&=P_{r}(d_0,X,tY^{-1}) \prod_{i=1}^{r} (1-t^{2}p^{-r-1+i})^2.
\end{align*}
This proves the assertion.

\end{proof}

\bigskip

 Now we consider a partial series of $\widetilde P_{r}(d_0,X,Y,t).$   For $d_0 \in {\bf Z}_p^*$, we put 
\begin{eqnarray*}
\lefteqn{Q_{r}(d_0,X,Y,t)} \\
 &=&  \sum_{B' \in \pi^{-i_p}\widetilde {\mathcal F}_{r,p}(d_0) \cap {\mathrm{Her}}_{r,*}({\mathcal O}_p)} {\widetilde G_p(\pi^{i_p}B',X,tY) \over \alpha_p(\pi^{i_p}B') } (tY^{-1})^{{\rm ord}(\det \pi^{i_p}B')}. 
 \end{eqnarray*}
  To consider the relation between $\widetilde P_{r}(d_0,X,Y,t) \ $ and $\ Q_{r}(d_0,X,Y,t),$ and to express $R_{m}(d_0,X,Y,t)$ in terms of 
$\widetilde P_{r}(d_0,X,Y,t),$ we provide some more preliminary results. 


Let $X$ be a variable. First suppose that $K_p$ is unramified over ${\bf Q}_p$ or $K_p={\bf Q}_p \oplus {\bf Q}_p.$ Put $\hat \xi_{p}= \sqrt{-1}$ or $1$ according as $K_p$ is unramified over ${\bf Q}_p$ or not. 
Let  $H_m=H_{m}(\cdot,X)$  be a function on  ${\mathrm{Her}}_m({\mathcal O}_p)^{\times}$  with values in ${\bf C}[X,X^{-1}]$ satisfying the following condition: 
$$H_{m}(1_{m-r} \bot pB,X)=\hat \xi_p^{(m-r){\rm ord}(\det (pB))}H_{r}(pB,\hat \xi_p^{m-r}X) \ {\rm  for\  any}  \ B \in {\mathrm{Her}}_{r}({\mathcal O}_p).$$
Let $d_0 \in {\bf Z}_p^*.$  Then we put
$$Q(d_0,H_{m},r,X,t)= \sum_{B \in p^{-1}{\mathcal F}_{r,p}(d_0) \cap {\mathrm{Her}}_{r} ({\mathcal O}_p)}{H_{m}(1_{m-r} \bot pB,X) \over \alpha_p(1_{m-r} \bot pB)}t^{{\rm ord}(\det (pB))}.$$
Next suppose that $K_p$ is ramified over ${\bf Q}_p.$  
Let  $H_{m}=H_m(\cdot,X)$  be a function on  ${\mathrm{Her}}_m({\mathcal O}_p)^{\times}$ with values in ${\bf C}[X,X^{-1}]$ satisfying the following condition: 

\bigskip

$H_{m}(\Theta_{m-r} \bot \pi^{i_p}B,X)=H_{r}(\pi^{i_p}B,X)$  for any  $B \in {\mathrm{Her}}_{r,*}({\mathcal O}_p)$ if $m-r$ is even. 

\bigskip

\noindent
Let $d_0 \in {\bf Z}_p^*$ and $m-r$ be even.   Then we put
$$Q(d_0,H_{m},r,X,t)=\sum_{B \in \pi^{-i_p}\widetilde {\mathcal F}_{r,p}(d_0) \cap {\mathrm{Her}}_{r,*} ({\mathcal O}_p)}{H_{m}(\Theta_{m-r} \bot \pi^{i_p}B,X)  \over \alpha_p(\Theta_{m-r} \bot \pi^{i_p}B)}t^{{\rm ord}(\det (\pi^{i_p}B))}.$$
Then we have  the following (cf.  [\cite{Kat3}, Proposition 4.2.4]). 

\bigskip

\noindent 
  {\bf Proposition 5.3.2.} {\it 
 
  \noindent 
  {\rm (1)} Suppose that $K_p$ is unramified over ${\bf Q}_p$ or $K_p={\bf Q}_p \oplus {\bf Q}_p.$ Then for any $d_0 \in {\bf Z}_p^*$ and a non-negative integer $r$ we have 
$$Q(d_0,H_{m},r,X,t)={Q(d_0,H_{r},r,\hat \xi_{p}^{m-r}X,\hat \xi_{p}^{m-r}t) \over \phi_{m-r}(\xi_p p^{-1})}.$$

  \noindent 
  {\rm (2)} Suppose that $K_p$ is ramified over ${\bf Q}_p.$  Then for any $d_0 \in {\bf Z}_p^*$ and a non-negative integer $r$ such that $m-r$ is even, we have 
$$Q(d_0,H_{m},r,X,t)={Q(d_0,H_{r},r,X,t) \over \phi_{(m-r)/2}(p^{-2})}.$$}

\bigskip

Now to apply Proposition 5.3.2  to the formal power series $R_{m}(d_0,X,Y,t) $ and $ Q_{r}(d_0, X,Y,t)$ we give  the following lemma.

\noindent
 {\bf Lemma 5.3.3.} {\it  
 Let $m$  be an integer.

  \noindent
  {\rm (1)}  Suppose that $K_p$ is unramified over ${\bf Q}_p$ or $K_p={\bf Q}_p \oplus {\bf Q}_p.$ Then for any integer such that $r \le m,$ and $B' \in {\mathrm{Her}}_{r}({\mathcal O}_p)$ we have   
  $$\widetilde G_p(1_{m-r} \bot  pB',X,t)=\widetilde G_p( pB', \hat \xi_{p}^{m-r}X,\hat \xi_{p}^{m-r}t).$$
  
 \noindent
{\rm (2)} Suppose that $K_p$ is ramified over ${\bf Q}_p.$  Then for any  non-negative integer $r$ such that $m-r$ is even, and  $B' \in {\mathrm{Her}}_{r,*}({\mathcal O}_p),$ we have  
  $$\widetilde G_p(\Theta_{m-r} \bot  \pi^{i_p}B',X,t)=\widetilde G_p(\pi^{i_p}B',X,t).$$

 }
  
  \begin{proof}    By Lemma 5.2.1 (1), we have
   $$G_p(1_{m-r} \bot  pB',X)=G_p(pB',\xi_p^{m-r}p^{m-r}X)$$   for $B' \in {\mathrm{Her}}_{r}({\mathcal O}_p).$ Hence  by Corollary to Lemma 5.2.2 we have 
   $$\widetilde F_p^{(0)}(1_{m-r} \bot  pB',X)=\hat \xi_{p}^{(m-r){\rm ord}(\det (pB'))}\widetilde F_p^{(0)}(pB',\hat \xi_{p}^{m-r}X)$$
   for $B' \in {\mathrm{Her}}_{r}({\mathcal O}_p).$ 
   Thus the assertion (1) follows from (3) of  Lemma 5.1.2.
   The assertion (2) can be proved in a similar way.
   \end{proof}

  \bigskip
  
 Let $R_{m}(d_0,X,Y,t)$ be the formal power series defined at the beginning of Section 5. We express $R_m(d_0,X,Y,t)$ in terms of   $Q_{r}(d_0,X,Y,t).$ 
 
\noindent    
 {\bf Theorem 5.3.4.} {\it  
Let $d_0 \in {\bf Z}_p^*.$

\noindent
{\rm (1)} Suppose that $K_p$ is unramified over ${\bf Q}_p.$ 
 Then  
\begin{align*}
R_{m}(d_0,X,Y,t)&=\sum_{r=0}^{m} {\prod_{i=0}^{r-1} (1-(-1)^m(-p)^{i}Y^2) \prod_{i=r}^{m-1}(1-(-1)^{m}(-p)^{-2m+i}Y^2t^2) \over \phi_{m-r}(-p^{-1})}\\
& \times  Q_{r}(d_0,\hat \xi_{p}^{m-r}X,p^{-m/2}Y,\hat \xi_{p}^{m-r}p^{-m/2}t).
\end{align*}

\noindent
{\rm (2)} Suppose that $K_p={\bf Q}_p \oplus {\bf Q}_p.$ 
 Then  
\begin{align*}
 R_{m}(d_0,X,Y,t)&=\sum_{r=0}^{m} {\prod_{i=0}^{r-1} (1-p^{i}Y^2) \prod_{i=r}^{m-1}(1-p^{-2m+i}Y^2t^2) \over \phi_{m-r}(p^{-1})}\\
& \times  Q_{r}(d_0,X,p^{-m/2}Y,p^{-m/2}t).
\end{align*}
Throughout {\rm (1)} and {\rm (2)}, we understand that $Q_0(d_0,X,Y,t)=1$.

\noindent
{\rm (3)} Suppose that $K_p$ is ramified over ${\bf Q}_p.$  Let $i_p=0$, or $1$ according as $p=2$ and $f_2=2$, or not as defined in Section 5.1. 

\begin{enumerate}
\item[{\rm (3.1)}] Let $m$ be odd. Then  
\begin{align*} 
 R_{m}(d_0,X,Y,t)&=\sum_{r=0}^{(m-1)/2} {\prod_{i=0}^{r-1} (1-p^{2i+1}Y^2) \prod_{i=r}^{(m-3)/2}(1-p^{-2m+2i+1}Y^2t^2) \over \phi_{(m-2r-1)/2}(p^{-2})}\\
& \times  (tY^{-1})^{(m-2r-1)i_p/2} Q_{2r+1}((-1)^{(m-2r-1)/2}d_0,X,p^{-m/2}Y,p^{-m/2}t).
\end{align*}

\item[{\rm (3.2)}] Let $m$ be even. Then  
\begin{align*} 
R_{m}(d_0,X,Y,t)&=\sum_{r=0}^{m/2} {\prod_{i=0}^{r-1} (1-p^{2i}Y^2) \prod_{i=r}^{(m-2)/2}(1-p^{-2m+2i}Y^2t^2) \over \phi_{(m-2r)/2}(p^{-2})}\\
& \times  (tY^{-1})^{(m-2r)i_p/2}Q_{2r}((-1)^{(m-2r)/2}d_0,X,p^{-m/2}Y,p^{-m/2}t).
\end{align*}
Here, for $u \in {\bf Z}_p^*$ we understand that $Q_0(u,X,Y,t)=1$ or $0$ according as $u \in N_{K_p/{\bf q}_p}({\mathcal O}_p^*)$ or not.

\end{enumerate}

}
  
\begin{proof}  First suppose that $K_p$ is unramified over ${\bf Q}_p$ or $K_p={\bf Q}_p \oplus {\bf Q}_p$. Let $B$ be an element of ${\widetilde{\mathrm{Her}}}_{r}({\mathcal O}_p).$ Then we note that $1_{m-r} \bot pB$ belongs to ${\widetilde{\mathcal F}}_{m,p}(d_0)$ if and only if $B \in p^{-1}{\widetilde{\mathcal F}}_{r,p}(d_0) \cap {\widetilde{\mathrm{Her}}}_{r}({\mathcal O}_p).$ Thus the assertions (1) and (2) follow from Lemmas 5.2.1, 5.2.3, and 5.3.3, and Proposition 5.3.2. 
 
 Next suppose that $K_p$ is ramified over ${\bf Q}_p.$ Let $B$ be an element of ${\widetilde{\mathrm{Her}}}_{r}({\mathcal O}_p).$ Let $m-r$ be even. Then we note that $\Theta_{m-r} \bot \pi^{i_p}B$ belongs to ${\widetilde{\mathcal F}}_{m,p}(d_0)$ if and only if $B \in \pi^{-i_p}{\widetilde{\mathcal F}}_{r,p}((-1)^{(m-r)/2}d_0) \cap {\mathrm{Her}}_{r,*}({\mathcal O}_p).$ Moreover we note that 
${\rm ord}(\det (\Theta_{m-r} \bot \pi^{i_p}B))=(m-r)i_p/2+{\rm ord}(\det (\pi^{i_p}B))$.  Thus the assertion (3) can be proved similarly to above.

 \end{proof}
  
 \bigskip
 
  Now to rewrite the above theorem, first we express $\widetilde P_{m}(d_0,X,Y,t)$ in terms of  $Q_{r}(d_0,X,Y,t).$ 
  
\noindent
{\bf Proposition 5.3.5.} {\it  
Let $d_0 \in {\bf Z}_p^*.$ 
   
    \noindent
{\rm (1)} Suppose that $K_p$ is unramified over ${\bf Q}_p$ or $K_p={\bf Q}_p \oplus {\bf Q}_p.$  Then 
  $$\widetilde P_{m}(d_0,\hat \xi_p^m X,Y,\hat \xi_p^m t) = \sum_{r=0}^{m}  {1 \over  \phi_{m-r}(\xi_p p^{-1})} Q_{r}(d_0,\hat \xi_p^r X,Y,\hat \xi_p^r t).$$
  
  \noindent
 {\rm (2)}  Suppose that $K_p$ is ramified over ${\bf Q}_p.$
 
 \noindent
 {\rm (2.1)} Let $m$ be odd. Then
\begin{align*}
 (tY^{-1})^{(1-m)i_p/2} \widetilde P_{m}((-1)^{(m-1)/2}d_0,X,Y,t) &= \sum_{r=0}^{(m-1)/2}  { 1 \over  \phi_{(m-2r-1)/2}(p^{-2})} \\
& \times (tY^{-1})^{-ri_p}Q_{2r+1}((-1)^rd_0,X,Y,t).
\end{align*}
  
 \noindent
 {\rm (2.2)} Let $m$ be even. Then
 \begin{align*}
(tY^{-1})^{-mi_p/2}\widetilde P_{m}((-1)^{m/2}d_0,X,Y,t) &= \sum_{r=0}^{m/2}  {1 \over  \phi_{(m-2r)/2}(p^{-2})}\\
& \times (tY^{-1})^{-ri_p} Q_{2r}((-1)^rd_0,X,Y,t).
\end{align*}}
\begin{proof}
 The assertion can be proved in the same argument as in the proof of Theorem 5.3.4. 
 \end{proof}

\noindent
{\bf Corollary.} {\it 
Let $d_0$ be an element of ${\bf Z}_{p}^*.$ 
 
  \noindent
  {\rm (1)} Suppose that $K_p$ is unramified over ${\bf Q}_p$ or $K_p={\bf Q}_p \oplus {\bf Q}_p.$  Then 
   $$Q_{r}(d_0,\hat \xi_p^r X,Y,\hat \xi_p^r t)=\sum_{m=0}^r {(-1)^m (\xi_p p)^{(m-m^2)/2} \over \phi_m(\xi_p p^{-1}) }\widetilde P_{r-m}(d_0,\hat \xi_p^{r-m} X,Y,\hat \xi_p^{r-m}t).$$
Here we understand that $\widetilde P_{0}(d_0,X,Y,t)=1$.

\noindent
{\rm (2)} Suppose that $K_p$ is ramified over ${\bf Q}_p.$
 Then  
 $$(tY^{-1})^{-ri_p} Q_{2r+1}((-1)^rd_0,X,Y,t)=\sum_{m=0}^r  {(-1)^m p^{m-m^2} \over \phi_m(p^{-2}) }(tY^{-1})^{(m-r)i_p}\widetilde P_{2r+1-2m}((-1)^{r-m}d_0,X,Y,t),$$ and
  $$(tY^{-1})^{-ri_p}Q_{2r}((-1)^rd_0,X,Y,t)=\sum_{m=0}^r  {(-1)^m p^{m-m^2} \over \phi_m(p^{-2}) }(tY^{-1})^{(m-r)i_p}\widetilde P_{2r-2m}((-1)^{r-m}d_0,X,Y,t).$$
Here, for $u \in {\bf Z}_p^*$ we understand that $\widetilde P_{0}(u,X,Y,t)=1$ or $0$ according as $u \in N_{K_p/{\bf Q}_p}({\mathcal O}_p^*)$ or not.
  }

\begin{proof}  We can prove the assertions by induction on $r$ (cf. [\cite{K-K4}, Corollary 5.1.2]).

\end{proof}

\bigskip

The following lemma follows from [\cite{I-K3}, Lemma 3.4].

\bigskip 

\noindent
{\bf Lemma 5.3.6.} {\it 
 Let $l$ be a positive integer. Then we have the following identity on the three variables $q,U$ and $Q:$  
	\begin{align*}
& \prod_{i=1}^l (1-U^{-1}Qq^{-i+1})U^l\\
 &=\sum_{m=0}^l {\phi_l(q^{-1}) \over \phi_{l-m}(q^{-1}) \phi_m(q^{-1})}\prod_{i=1}^{l-m} (1-Qq^{-i+1})\prod_{i=1}^m (1-Uq^{i-1}) (-1)^mq^{(m-m^2)/2}.
\end{align*}
  }

\noindent 
  {\bf Theorem 5.3.7.} {\it 
  Let the notation be as in Theorem 5.3.5. 
  
  \noindent
{\rm (1)} Suppose that $K_p$ is unramified over ${\bf Q}_p$ or $K_p={\bf Q}_p \oplus {\bf Q}_p.$  Then 
\begin{align*} 
&R_{m}(d_0,X,Y,t) = \sum_{l=0}^{m} ((p^l\xi_pY^2)^{m-l}\widetilde P_{l}(d_0,\hat \xi_p^{m-l}X,p^{-m/2}Y,\hat \xi_p^{m-l}p^{-m/2}t) \\
& \times  {\prod_{i=1}^{m-l}(1- (\xi_p p)^{-l-m-i}t^2)\prod_{i=0}^{l-1}(1- \xi_p^m(\xi_p p)^{i}Y^2) \over \phi_{m-l}(\xi_p p^{-1})}. 
\end{align*}

\noindent
{\rm (2)} Suppose that $K_p$ is ramified over ${\bf Q}_p.$

\noindent
{\rm (2.1)} Let $m$ be odd.  Then   
\begin{align*}
&R_{m}(d_0,X,Y,t)=\sum_{l=0}^{(m-1)/2} (tY^{-1})^{(m-2l-1)i_p/2} \widetilde P_{2l+1}((-1)^{(m-2l-1)/2}d_0,X,p^{-m/2}Y,p^{-m/2}t)   \\
& \times {(p^{2l+1}Y^2)^{(m-2l-1)/2} \prod_{i=0}^{l-1} (1-p^{2i+1}Y^2)\prod_{i=1}^{(m-2l-1)/2}(1-p^{-2l-m-2i-1}t^2) \over \phi_{(m-2l-1)/2}(p^{-2})}.
\end{align*}
\noindent
{\rm (2.2)} Let $m$ be even.  Then   
\begin{align*}
& R_{m}(d_0,X,Y,t)=\sum_{l=0}^{m/2}  (tY^{-1})^{(m-2l)i_p/2}\widetilde P_{2l}((-1)^{(m-2l)/2}d_0,X,p^{-m/2}Y,p^{-m/2}t)   \\
& \times {(p^{2l}Y^2)^{(m-2l)/2} \prod_{i=0}^{l-1} (1-p^{2i}Y^2)\prod_{i=1}^{(m-2l)/2}(1-p^{-2l-m-2i}t^2) \over \phi_{(m-2l)/2}(p^{-2})}.
\end{align*}

}

\begin{proof} (1) By Theorem 5.3.4 and Corollary to Proposition 5.3.5, we have
\begin{align*}
& R_{m}(d_0,X,Y,t)\\
&=\sum_{r=0}^{m} {\prod_{i=0}^{r-1} (1-\xi_p^m(\xi_pp)^{i}Y^2) \prod_{i=0}^{m-r-1}(1-(\xi_p p)^{-m+i+r}p^{-m}Y^2t^2) \over \phi_{m-r}((\xi_pp)^{-1})}\\
& \times \sum_{j=0}^r  {(-1)^j (\xi_p p)^{(j-j^2)/2} \over \phi_j((\xi_p p)^{-1})} \widetilde P_{r-j}(d_0,\hat \xi_p^{m-r+j}X,p^{-m/2}Y,\hat \xi_p^{m-r+j}p^{-m/2}t)\\
&=\sum_{l=0}^{m} \widetilde P_{l}(d_0,\hat \xi_p^{m-l}X,p^{-m/2}Y,\hat \xi_p^{m-l}p^{-m/2}t)\\
&\times \sum_{j=0}^{m-l} (-1)^j (\xi_p p)^{(j-j^2)/2}{\prod_{i=0}^{l+j-1} (1-\xi_p^m(\xi_pp)^{i}Y^2) \prod_{i=0}^{m-l-j-1}(1-(\xi_p p)^{-m+i+l+j}p^{-m}Y^2t^2) \over \phi_j(\xi_p p^{-1})\phi_{m-j-l}(\xi_p p^{-1})}.
\end{align*}
Then  the assertion (1) follows from Lemma 5.3.6. 

 (2) Let $m$ be odd. Then, again by Theorem 5.3.4 and Corollary to Proposition 5.3.5, 
\begin{align*}
&R_m(d_0,X,Y,t)\\
&=\sum_{r=0}^{(m-1)/2} {\prod_{i=0}^{r-1}(1-p^{2i+1}Y^2) \prod_{i=0}^{(m-1)/2-r-1}(1-p^{-2m+2i+2r+1}Y^2t^2) \over \phi_{(m-2r-1)/2}(p^{-2})}\\
& \times (tY^{-1})^{(m-1)i_p/2} \sum_{j=0}^r {(-1)^jp^{j-j^2} \over \phi_j(p^{-2})} (tY^{-1})^{(j-r)i_p} \\
&\times \widetilde P_{2r+1-2j}((-1)^{(m-1-2r+2j)/2}d_0,X,p^{-m/2}Y,p^{-m/2}t)\\
&=(tY^{-1})^{(m-1)i_p/2} \sum_{l=0}^{(m-1)/2} (tY^{-1})^{-l i_p}\widetilde P_{2l+1}((-1)^{(m-1-2l)/2}d_0,X,Y^{-m/2}Y,p^{-m/2}t)\\
&\times \sum_{j=0}^{(m-1)/2-l} (-1)^jp^{j-j^2} {\prod_{i=0}^{l+j-1}(1-p^{2i+1}Y^2)\prod_{i=0}^{(m-1)/2-l-j-1}(1-p^{-2m+2i+2l+2j+1}Y^2t^2)
\over \phi_j(p^{-2})\phi_{(m-1)/2-j-l}(p^{-2})}.
\end{align*}
Hence the assertion (2.1) follows from Lemma 5.3.6. The assertion (2.2) can be proved  in the same manner as above.
\end{proof}
By Proposition 5.3.1 we obtain:

\noindent
{\bf Corollary.} {\it
 \noindent
{\rm (1)} Suppose that $K_p$ is unramified over ${\bf Q}_p$ or $K_p={\bf Q}_p \oplus {\bf Q}_p.$  Then 
\begin{align*}
&R_{m}(d_0,X,Y,t) =\prod_{i=1}^m (1-p^{-2m}(\xi_p p)^{i-1}t^2)\\
& \times \sum_{l=0}^{m} (p^l\xi_pY^2)^{m-l} P_{l}(d_0,\hat \xi_p^{m-l}X,\hat \xi_p^{m-l}tY^{-1})  {\prod_{i=1}^{l}(1-\xi_p  (\xi_p p)^{-l-m+i-1}t^2)\prod_{i=0}^{l-1}(1- \xi_p^m(\xi_p p)^{i}Y^2) \over \phi_{m-l}(\xi_p p^{-1})}. 
\end{align*}
Here we understand that $P_{0}(d_0,X,t)=1$.

\noindent
{\rm (2)} Suppose that $K_p$ is ramified over ${\bf Q}_p.$

\noindent
{\rm (2.1)} Let $m$ be odd.  Then   
\begin{align*}
& R_{m}(d_0,X,Y,t)=\prod_{i=1}^{(m+1)/2}(1-p^{-2m+2i-2}t^2)\\
&\times  \sum_{l=0}^{(m-1)/2} (tY^{-1})^{(m-2l-1)i_p/2} P_{2l+1}((-1)^{(m-2l-1)/2}d_0,X,tY^{-1})   \\
& \times {(p^{2l+1}Y^2)^{(m-2l-1)/2} \prod_{i=0}^{l-1} (1-p^{2i+1}Y^2)\prod_{i=1}^{l}(1-p^{-2l-2+2i-m}t^2) \over \phi_{(m-2l-1)/2}(p^{-2})}.
\end{align*}
\noindent
{\rm (2.2)} Let $m$ be even.  Then   
\begin{align*}
& R_{m}(d_0,X,Y,t)=\prod_{i=1}^{m/2}(1-p^{-2m+2i-2}t^2) \\
&\times \sum_{l=0}^{m/2} (tY^{-1})^{(m-2l)i_p/2} P_{2l}((-1)^{(m-2l)/2}d_0,X,tY^{-1})   \\
& \times {(p^{2l}Y^2)^{(m-2l)/2} \prod_{i=0}^{l-1} (1-p^{2i}Y^2)\prod_{i=1}^{l}(1-p^{-2l-1+2i-m}t^2) \over \phi_{(m-2l)/2}(p^{-2})}.
\end{align*}
Here, for $u \in {\bf Z}_p^*$ we understand that $P_{0}(u,X,t)=1$ or $0$ according as $u \in N_{K_p/{\bf Q}_p}({\mathcal O}_p^*)$ or not. 

}

 \subsection{Explicit formulas of formal power series of Koecher-Maass type}
 
 { }
 \noindent
 { }
 
 \bigskip
 
 In this section we review explicit formulas for $P_m(d_0,X,t).$  

\noindent
 {\bf Theorem 5.4.1.} {\it 
{\rm [\cite{Kat3}, Theorem 4.3.1]} Let $m$ be even, and $d_0 \in {\bf Z}_{p}^*.$ 
 
\noindent 
 {\rm (1)} Suppose that $K_p$ is unramified over ${\bf Q}_p.$   Then 
 $$P_m(d_0,X,t)={ 1 \over \phi_{m}(-p^{-1})\prod_{i=1}^{m} (1-t (-p)^{-i}X)(1+t(-p)^{-i}X^{-1})  }.$$

\noindent  
  {\rm (2)} Suppose that $K_p={\bf Q}_p \oplus {\bf Q}_p.$ Then 
  $$P_m(d_0,X,t)={ 1 \over \phi_{m}(p^{-1})\prod_{i=1}^{m} (1-t p^{-i}X)(1-tp^{-i}X^{-1})  }.$$

\noindent 
 {\rm (3)} Suppose that $K_p$ is ramified over ${\bf Q}_p.$ Let $\chi_{K_p}$ be the character of ${\bf Q}_p^:$ defined by
$\chi_{K_p}(a)=(-D,a)$ for $a \in {\bf Q}_p^*.$  Then 
\begin{align*}
& P_m(d_0,X,t)={ t^{mi_p/2} \over 2\phi_{m/2}(p^{-2})}\\
&\times \left\{{1 \over \prod_{i=1}^{m/2} (1-t p^{-2i+1}X)(1-tp^{-2i}X^{-1})} + {\chi_{K_p}((-1)^{m/2}d_0)  \over \prod_{i=1}^{m/2} (1-t p^{-2i}X)(1-tp^{-2i+1}X^{-1})}\right \}.
\end{align*}

  }
   
\noindent
  {\bf Theorem 5.4.2.} {\it  
 {\rm [\cite{Kat3}, Theorem 4.3.2]} Let $m$ be odd, and $d_0 \in {\bf Z}_{p}^*.$ 
 
\noindent 
 {\rm (1)} Suppose that $K_p$ is unramified over ${\bf Q}_p.$   Then 
 $$P_m(d_0,X,t)={ 1 \over \phi_{m}(-p^{-1})\prod_{i=1}^{m} (1+t (-p)^{-i}X)(1+t(-p)^{-i}X^{-1})  }.$$
 
\noindent  
  {\rm (2)} Suppose that $K_p={\bf Q}_p \oplus {\bf Q}_p.$ Then 
  $$P_m(d_0,X,t)={ 1 \over \phi_{m}(p^{-1})\prod_{i=1}^{m} (1-t p^{-i}X)(1-tp^{-i}X^{-1})  }.$$

\noindent 
 {\rm (3)} Suppose that $K_p$ is ramified over ${\bf Q}_p.$   Then 
 $$P_m(d_0,X,t)={ t^{(m+1)i_p/2+\delta_{2p}} \over 2\phi_{(m-1)/2}(p^{-2})\prod_{i=1}^{(m+1)/2} (1-t p^{-2i+1}X)(1-tp^{-2i+1}X^{-1})  }.$$
   } 

\subsection{Explicit formulas of formal power series of Rankin-Selberg type}
 
 \noindent
  { }
  
  \bigskip
  
  We give an  explicit formula for $H_m(d,X,Y,t).$  First we remark the following.

\noindent
 {\bf Proposition 5.5.1.} {\it 
Let $d \in {\bf Z}_p^{\times}.$ 
 Then we have 
$$\lambda_{m,p}^*(d,X,Y)=u_p\lambda_{m,p}(d,X,Y).$$ }
   
 \begin{proof}
 This can be proved in the same way as [\cite{Kat3}, Proposition 4.3.7]
 \end{proof}

It is well known that $\#({\bf Z}_p^*/N_{K_p/{\bf Q}_p}({\mathcal O}_p^*))=2$ if $K_p/{\bf Q}_p$ is ramified. Hence we can take a complete set ${\mathcal N}_p$ of representatives of ${\bf Z}_p^*/N_{K_p/{\bf Q}_p}({\mathcal O}_p^*)$ so that ${\mathcal N}_p=\{1,\xi_0\}$ with $\chi_{K_p}(\xi_0)=-1.$

 \bigskip

\noindent
 {\bf Theorem 5.5.2.} {\it 
Let $m=2n$ be even, and  $d_0 \in {\bf Z}_p^*.$ 

\noindent
  {\rm (1)}  Suppose that $K_p$ is unramified over ${\bf Q}_p.$ Then 
\begin{align*}
 &H_{2n}(d_0,X,Y,t)= {\prod_{i=1}^{2n} (1-p^{-4n}(-p)^{i-1}t^2) \over \phi_{2n}(-p^{-1})} \\
 & \times  { 1 \over \prod_{i=1}^{2n} (1+(-p)^{-2n+i-1}XYt)(1-(-p)^{-2n+i-1}XY^{-1}t)}\\
 & \times {1 \over \prod_{i=1}^{2n} (1-(-p)^{-2n+i-1}X^{-1}Yt)(1+(-p)^{-2n+i-1}X^{-1}Y^{-1}t)}.
\end{align*}
   {\rm (2)}  Suppose that $K_p={\bf Q}_p \oplus {\bf Q}_p.$ Then 
\begin{align*}
 & H_{2n}(d_0,X,Y,t)= {\prod_{i=1}^{2n} (1-p^{-4n}p^{i-1}t^2) \over \phi_{2n}(p^{-1})} \\
& \times  { 1 \over \prod_{i=1}^{2n} (1-p^{-2n+i-1}XYt)(1-p^{-2n+i-1}XY^{-1}t)}\\
 & \times {1 \over \prod_{i=1}^{2n} (1-p^{-2n+i-1}X^{-1}Yt)(1-p^{-2n+i-1}X^{-1}Y^{-1}t)}.
\end{align*}
   {\rm (3)}  Suppose that $K_p$ is ramified over ${\bf Q}_p.$ 
For $l=0,1$ put
$$H_{2n}^{(l)}(X,Y,t)=\sum_{d \in {\mathcal N}_p} \chi_{K_p}((-1)^nd)^l H_{2n}(d,X,Y,t).$$
Then we have
$$H_{2n}(d_0,X,Y,t)={1 \over 2} (H_{2n}^{(0)}(X,Y,t)+\chi_{K_p}((-1)^nd_0)H_{2n}^{(1)}(X,Y,t)),$$
and 
\begin{align*}
 &H_{2n}^{(0)}(X,Y,t)= t^{ni_p}{\prod_{i=1}^{n} (1-p^{-4n}p^{2i-2}t^2) \over \phi_n(p^{-2})} \\
 & \times {  1 \over \prod_{i=1}^{n} (1-p^{-2n+2i-1}XYt)(1-p^{-2n+2i-1}X^{-1}Y^{-1}t)}\\
 &\times  { 1 \over \prod_{i=1}^{n} (1-p^{-2n+2i-2}X^{-1}Yt)(1-p^{-2n+2i-2}XY^{-1}t)},
\end{align*}
 and
\begin{align*}
 &H_{2n}^{(1)}(X,Y,t)=t^{ni_p}{\prod_{i=1}^{n} (1-p^{-4n}p^{2i-2}t^2) \over \phi_n(p^{-2})} \\
 &\times { 1 \over \prod_{i=1}^{n} (1-p^{-2n+2i-1}X^{-1}Yt)(1-p^{-2n+2i-1}XY^{-1}t)} \\
&\times  { 1 \over \prod_{i=1}^{n} (1-p^{-2n+2i-2}XYt)(1-p^{-2n+2i-2}X^{-1}Y^{-1}t)} \}.
\end{align*}

 }
 
\begin{proof} First we prove (1). 
By Theorems 5.4.1and 5.4.2, we have  
$$P_l(d_0,\hat \xi_p^{m-l}X,\hat \xi_p^{m-l}X)=P_l(d_0,X,t)$$
if $l$ is even, and
$$P_l(d_0,\hat \xi_p^{m-l}X,\hat \xi_p^{m-l}X)={1 \over \phi_m(-p^{-1}) \prod_{i=1}^l (1-t(-p)^{-i}X)(1+t(-p)^{-i}X^{-1})}$$
if $l$ is odd. Hence, by Corollary to Theorem 5.3.7, $R_{2n}(d_0,X,Y,t)$ can be expressed as
\begin{align*}
&R_{2n}(d_0,X,Y,t)\\
&= {\prod_{i=1}^{2n}(1-p^{-4n}(-p)^{i-1}t^2) S(X,Y,t) \over \phi_{2n}(-p)\prod_{i=1}^{2n}(1-t(-p)^{-2n+i-1}XY^{-1})(1+t(-p)^{-2n+i-1}X^{-1}Y^{-1})},
\end{align*}
where $S(X,Y,t)$ is a polynomial in $t$ of degree at most $4n.$ 
Then by Theorem 5.2.8, we have
\begin{align*}
&H_{2n}(d_0,X,Y,t)\\
&={\prod_{i=1}^{2n}(1-p^{-4n}(-p)^{i-1}t^2) S(X,Y,t) \over \phi_{2n}(-p)\prod_{i=1}^{2n}(1-t(-p)^{-2n+i-1}XY^{-1})(1+t(-p)^{-2n+i-1}X^{-1}Y^{-1})}\\
& \times {1 \over \prod_{i=1}^{2n}(1-t^2p^{-4n+2i-2}X^2Y^2)(1-t^2p^{-4n+2i-2}X^{-2}Y^{2})}.
\end{align*}
Recall that we have the following functional equation
$$H_{2n}(d_0,X,Y^{-1},t)=H_{2n}(d_0,X,-Y,t).$$ 
Hence the reduced denominator of the rational function $H_{2n}(d_0,X,Y^{-1},t)$ in $t$ is at most
\begin{align*}
&\prod_{i=1}^{2n}\{(1-t(-p)^{-2n+i-1}XY^{-1})(1+t(-p)^{-2n+i-1}X^{-1}Y^{-1})\\
& \times (1+t(-p)^{-2n+i-1}XY)(1-t(-p)^{-2n+i-1}X^{-1}Y)\},
\end{align*}
and therefore we have
\begin{align*}
&H_{2n}(d_0,X,Y,t)={c \prod_{i=1}^{2n}(1-(-p)^{-2n-i}t^2) \over  \phi_{2n}(-p)}\\
& \times {1 \over \prod_{i=1}^{2n}(1-t(-p)^{-2n+i}XY^{-1})(1+t(-p)^{-2n+i}X^{-1}Y^{-1})}\\
&\times {1 \over \prod_{i=1}^{2n}(1+t(-p)^{-2n+i-1}XY)(1-t(-p)^{-2n+i-1}X^{-1}Y)}
\end{align*}
with some constant $c.$ We easily see that we have $c=1.$ This proves the assertion (1).
Similarly the assertions (2) and (3) can be proved.
\end{proof}

\bigskip

\bigskip

Similarly to Theorem 5.5.2, we have

\bigskip
 
\noindent  
 {\bf Theorem 5.5.3.} {\it Let $m=2n+1$ be odd, and  $d_0 \in {\bf Z}_p^*.$ 

\noindent
  {\rm (1)}  Suppose that $K_p$ is unramified over ${\bf Q}_p.$ Then 
\begin{align*}
 &H_{2n+1}(d_0,X,Y,t)= {\prod_{i=1}^{2n+1} (1-p^{-4n-2}(-p)^{i-1}t^2) \over \phi_{2n+1}(-p^{-1})} \\
 &\times  { 1 \over \prod_{i=1}^{2n+1} (1+(-p)^{-2n+i-2}XYt)(1+(-p)^{-2n+i-2}XY^{-1}t)}\\
 & \times {1 \over \prod_{i=1}^{2n} (1+(-p)^{-2n+i-2}X^{-1}Yt)(1+(-p)^{-2n+i-2}X^{-1}Y^{-1}t)}.
\end{align*} 
 
 \noindent
  {\rm (2)}  Suppose that $K_p={\bf Q}_p \oplus {\bf Q}_p.$ Then 
\begin{align*}
 &H_{2n+1}(d_0,X,Y,t)={\prod_{i=1}^{2n+1} (1-p^{-4n-2}p^{i-1}t^2) \over \phi_{2n+1}(p^{-1})} \\
 & \times  { 1 \over \prod_{i=1}^{2n+1} (1-p^{-2n+i-2}XYt)(1-p^{-2n+i-2}XY^{-1}t)}\\
 & \times {1 \over \prod_{i=1}^{2n+1} (1-p^{-2n+i-2}X^{-1}Yt)(1-p^{-2n+i-2}X^{-1}Y^{-1}t)}.
\end{align*}

  \noindent
  {\rm (3)}  Suppose that $K_p$ is ramified over ${\bf Q}_p.$  Then 
\begin{align*}
 &H_{2n+1}(d_0,X,Y,t)= t^{(n+1)i_p+\delta_{2p}}{\prod_{i=1}^{n+1} (1-p^{-4n-2}p^{2i-2}t^2) \over 2 \phi_n(p^{-2})} \\
 & \times  { 1 \over \prod_{i=1}^{n+1} (1-p^{-2n+2i-3}XYt)(1-p^{-2n+2i-3}X^{-1}Y^{-1}t)}\\
 & \times {1 \over  (1-p^{-2n+2i-3}X^{-1}Yt)(1-p^{-2n+2i-3}XY^{-1}t)   }.
\end{align*}
  }

\bigskip
By using the same argument as in the proof of [\cite{Kat3},Theorem 4.3.6 and its corollary] we obtain the following:

\bigskip

\noindent
{\bf Theorem 5.5.4.} {\it 
Let $d_0 \in {\bf Z}_p^*.$ 

\noindent
{\rm (1)}  Suppose that  $K_p$ is unramified  over ${\bf Q}_p$ or that $K_p={\bf Q}_p \oplus {\bf Q}_p.$ Then
$$\hat H_{m}(d_0,X,Y,t)=H_{m}(d_0,X,Y,t)$$
for any $m >0.$

\noindent
{\rm (2)}  Suppose that  $K_p$ is ramified  over ${\bf Q}_p.$  \\
{\rm (2.1)} For $l=0,1$ put
$$\hat H_{2n}^{(l)}(X,Y,t)=\sum_{d \in {\mathcal N}_p} \chi_{K_p}((-1)^nd)^l \hat H_m(d,X,Y,t).$$
Then we have
$$\hat H_{2n}(d_0,X,Y,t)={1 \over 2} (\hat H_{2n}^{(0)}(X,Y,t)+\chi_{K_p}((-1)^nd_0)\hat H_{2n}^{(1)}(X,Y,t)),$$
and 
 $$\hat H_{2n}^{(0)}(X,Y,t)=H_{2n}^{(0)}(X,Y,t),$$
and
$$\hat H_{2n}^{(1)}(X,Y,t)=H_{2n}^{(1)}(X,Y,\chi_{K_p}(p) t).$$
{\rm (2.2)} We have
$$\hat H_{2n+1}(d_0,X,Y,t)=H_{2n+1}(d_0,X,Y,t)$$

}
\section{Proof of the main theorem}

\noindent
{\bf Theorem 6.1.} {\it 
Let $k$ and $n$ be positive integers. Let $f$ be a primitive form in ${\mathfrak S}_{2k+1}(\varGamma_0(D),\chi).$ For a subset $Q$ of $Q_D$ and a Dirichlet character $\eta=\chi^{i-1}$ with a positive integer $i$ put
\begin{align*}
&M(s,f,{\rm Ad},\eta,\chi_Q)\\
&=\{\prod_{p \not\in Q} (1-\alpha_p^2 \chi(p)^i \chi_Q(p) p^{-s})(1-\alpha_p^{-2} \chi(p)^i \chi_Q(p) p^{-s})(1-\chi^{i-1}(p)\chi_Q(p)p^{-s})^2 \\
&\times \prod_{p \in Q} (1-\alpha_p^2 \chi_Q'(p)\chi^{i-1}(p)p^{-s})(1-\alpha_p^{-2} \chi_Q'(p)\chi^{i-1}(p) p^{-s})(1-\chi_Q'(p)\chi(p)^i p^{-s})^2 \}^{-1},
\end{align*}
where for $\psi=\chi_Q$ or $\psi=\chi_Q'$ we make the convention $\psi(p)\chi^j(p)=\psi(p)$ or $0$ according as $j$ is even or odd. Then, we  have
\begin{align*}
& R(s,I_{2n}(f))=D^{ns+n^2-n/2-1/2}2^{-2n+1}\\
& \times \prod_{i=2}^{2n} \widetilde {\Lambda}(i,\chi^i)\prod_{i=0}^{2n-1}L(2s-4k-i,\chi^i)^{-1}\\
& \times \sum_{Q \subset Q_D } \chi_Q((-1)^n)\prod_{i=1}^{2n}M(s-2k-2n+i,f,{\rm Ad},\chi^{i-1},\chi_Q) .
\end{align*}

}

\begin{proof} The assertion can be proved by using Theorems 4.1, 5.5.2 and 5.5.4
similarly to [\cite{Kat3}, Theorem 2.3].
 
\end{proof}

\noindent
{\bf Theorem 6.2.} {\it 
Let $k$ and $n$ be positive integers. Given a primitive form $f \in {\mathfrak S}_{2k}(SL_2({\bf Z})).$ Then, we  have
\begin{align*}
&R(s,I_{2n+1}(f))=D^{ns+n^2+3n/2+1/2}2^{-2n}\\
&\times \prod_{i=2}^{2n+1} \widetilde {\Lambda}(i,\chi^i)\prod_{i=0}^{2n}L(2s-4k-i+2,\chi^i)^{-1}\\
&\times \prod_{i=1}^{2n+1}L(s-2k-2n+i,f,{\rm Ad},\chi^{i-1})L(s-2k-2n+i,\chi^{i-1}).
\end{align*}
}

\begin{proof} The assertion follows directly from Theorems 4.1 and 5.5.3.
\end{proof}

\noindent
{\bf Lemma 6.3.} {\it 
Let $f$ be a primitive form in ${\mathfrak S}_{2k+1}(\varGamma_0(D),\chi).$ Suppose that 
$f_Q=f$ for $Q \subset Q_D.$ Then for a positive integer $i$ we have
$$M(s,f,{\rm Ad},\chi^{i-1},\chi_Q) ={L}(s,f,{\rm Ad},\chi^{i-1}) {L}(s,\chi^{i-1}).$$}

\begin{proof} For a prime number $p$ let $M_p(s)$ and $L_p(s)$ be the $p$-Euler factor of $M(s,f,{\rm Ad},\chi^{i-1},\chi_Q)$ and ${L}(s,f,{\rm Ad},\chi^{i-1}) {L}(s,\chi^{i-1}),$
 respectively. We have  $M_p(s)=L_p(s)$ if $p \not\in Q$ and $\chi_Q(p)=1.$ 
 By the assumption  we have 
 $$\chi_Q(p)c_f(p)=c_f(p).$$
 Since $f$ is a primitive form, we have  $c_f(p) \not=0$ for $p | D.$ Hence we have
 $M_p(s)=L_p(s)$ if $p \not\in Q$ and $p |D.$ Suppose $p \nmid D$ and $\chi_Q(p)=-1.$ Then $c_f(p)=0$ and hence $\alpha_p +\chi(p)\alpha_p^{-1}=0.$
 Then by a simple computation we have
 $$M_p(s)=(1-p^{-2s})^{-2}.$$
 Similarly we have
 $$ L_p(s)=(1-p^{-2s})^{-2}.$$
Suppose that $p \in Q.$ 
 Then $|\alpha_p|= |c_f(p)|=1,$  and $\chi_Q'(p)\overline {c_f(p)}=c_f(p).$ Hence
 $\alpha_p$ is a real number or a purely imaginary number according as $\chi_Q'(p)=1$ or $-1.$ Hence
 $\chi_Q'(p)\alpha_p^2=\chi_Q'(p)\alpha_p^{-2}=1,$ and
 $$M_p(s)=L_p(s).$$ 
This completes the assertion.
 \end{proof}
 
 \noindent
 {\bf Proposition 6.4.} {\it 
  {\rm (1)} Let $f$ be a primitive form in ${\mathfrak S}_{2k+1}(\varGamma_0(D),\chi),$ and $Q$ be a subset of $Q_D.$   Then for a positive integer $i \ge 2$ the Euler product $M(s+i-1,f,{\rm Ad},\chi^{i-1},\chi_Q)$ is holomorphic at $s=1.$ Moreover $M(s,f,{\rm Ad},1,\chi_Q) $ has a  non-zero residue at $s=1$ if and only if $f=f_Q.$ In this case the residue of $M(s,f,{\rm Ad},1,\chi_Q)$ at $s=1$ is ${L}(1,f,{\rm Ad}).$
  
  {\rm (2)} Let $f$ be a primitive form in ${\mathfrak S}_{2k}(SL_2({\bf Z}))$ and $\chi$ be a primitive quadratic odd character. Then for a positive integer $i \ge 2$ the Euler product ${L}(s+i-1,f,{\rm Ad},\chi^{i-1}) {L}(s+i-1,\chi^{i-1})$ is holomorphic at $s=1,$ and ${L}(s,f,{\rm Ad},1) {L}(s,1)$ has a simple pole at $s=1$ with the residue ${L}(1,f,{\rm Ad}).$}
   
\begin{proof} (1) Clearly $M(s+i-1,f,{\rm Ad},\chi^{i-1},\chi_Q)$ is holomorphic at $s=1$ if $ i \ge 2.$
To prove the latter half of the assertion, let 
$R(s, f_Q \otimes f_{\rho})$ be the tensor product $L$-function of $f_Q$ and $f_{\rho},$ where
$$f_{\rho}(z)=\sum_{e=1}^{\infty} \overline {c_f(e)}{\bf e}(ez).$$  
We note that
$\overline {c_f(e)}=\chi(e)c_f(n)$ and $c_{f_Q}(e)=\chi_Q(e)c_f(n)$  if $(e,D)=1.$ 
Hence we have
$$M(s,f,{\rm Ad},1,\chi_Q)=R(s, f_Q \otimes f_{\rho}) \times \prod_{p | D} {M_p(s,f,{\rm Ad},1,\chi_Q) \over R_p(s, f_Q \otimes f_{\rho})},$$
where $M_p(s,f,{\rm Ad},1,\chi_Q)$ and $R_p(s, f_Q \otimes f_{\rho})$ are the $p$-Euler factors of 
$M(s,f,{\rm Ad},1,\chi_Q)$ and  $R(s, f_Q \otimes f_{\rho}),$ respectively.  We note 
$\prod_{p | D} {M_p(s,f,{\rm Ad},1,\chi_Q) \over R_p(s, f_Q \otimes f_{\rho})}$ is holomorphic and nonzero at $s=1.$
Hence we have
$${\rm Res}_{s=1}M(s,f,{\rm Ad},1,\chi_Q) =c(f_Q,f)$$
with $c$ a nonzero complex numbers (cf. [\cite {Sh0}, p. 788] and [\cite {Sh3}, p. 831]).
Hence $M(s,f,{\rm Ad},1,\chi_Q) $ has a  non-zero residue at $s=1$ if and only if
 $(f,f_Q) \not=0.$ Since $f$ and $f_Q$ are primitive forms, this is equivalent to say that $f=f_Q.$ In this case, 
we have 
 $$M(s,f,{\rm Ad},1,\chi_Q)=L(s,f,{\rm Ad})\zeta(s),$$
 and hence the last assertion holds.
 
 (2) The assertion can easily be proved.

\end{proof}

\bigskip  
{\bf Proof of Theorem 2.1.}\\
 (1) By Theorem 6.1 and Lemma 6.3, we have
\begin{align*}
 &R(s,I_m(f))=D^{ns+n^2-n/2-1/2}2^{-2n+1}\prod_{i=1}^{2n} \widetilde {\Lambda}(i,\chi^i)\prod_{i=0}^{2n-1}L(2s-4k-i,\chi^i)^{-1}\\
& \times \{\eta_m(f) \prod_{i=1}^{2n}L(s-2k-2n+i,f,{\rm Ad},\chi^{i-1})L(s-2k-2n+i,\chi^{i-1})\\
&+ \sum_{Q \in Q_D \atop f_Q \not= f} \chi_Q((-1)^n)\prod_{i=1}^{2n}M(s-2k-2n+i,f,{\rm Ad},\chi^{i-1},\chi_Q) \}.
\end{align*}
By (1) of  Lemma 6.4, the term
$$\prod_{i=0}^{2n-1}L(2s-4k-i,\chi^i)^{-1} \prod_{i=1}^{2n}M(2s-2k+i,f,{\rm Ad},\chi^{i-1},\chi_Q)$$
is holomorphic at $s=2k+2n$ if $f_Q \not= f.$ On the other hand, the term
$$\prod_{i=0}^{2n-1}L(2s-4k-i,\chi^i)^{-1}\prod_{i=1}^{2n}L(s-2k-2n+i,f,{\rm Ad},\chi^{i-1})L(s-2k-2n+i,\chi^{i-1})$$
 has a simple pole at $s=2k+2n$ with the residue
 $$\prod_{i=0}^{2n-1}L(4n-i,\chi^i)^{-1}\prod_{i=1}^{2n}L(i,f,{\rm Ad},\chi^{i-1}) \prod_{i=2}^{2n}L(i,\chi^{i-1}).$$
 Hence $R(s,I_m(f))$ has a simple at $s=2k+2n$ with the residue 
\begin{align*}
&D^{n(2k+2n)+n^2-n/2-1/2}2^{-2n+1}\\
 & \times \eta_m(f)\prod_{i=2}^{2n} \widetilde {\Lambda}(i,\chi^i)\prod_{i=0}^{2n-1}L(4n-i,\chi^i)^{-1}\prod_{i=1}^{2n}L(i,f,{\rm Ad},\chi^{i-1}) \prod_{i=2}^{2n}L(i,\chi^{i-1}).
\end{align*}
  Thus the assertion can be proved by comparing the above result with Proposition 3.1. 
 
 (2) The assertion holds if $m=1.$ In the case $m \ge 3,$ the assertion can be proved by Theorem 6.2, (2) of Lemma 6.4,  and Proposition 3.1 in the same manner as above.
 
 \bigskip
{\bf Acknowledgement}
The author was partly supported by the JSPS KAKENHI Grant Numbers 24540005, 25247001 and 23224001.
The author  thanks T. Ikeda for useful comments. The author also thanks the referee for pointing out many errors in the
original  version of our paper.

\end{document}